\documentclass[11pt,a4paper]{article}

\usepackage{soul}
\usepackage[normalem]{ulem}

\usepackage{epsf,epsfig,amsfonts,amsgen,amsmath,amstext,amsbsy,amsopn,lineno}
\usepackage{amsmath,amsthm,amsbsy}
\usepackage{color}
\usepackage{cite}
\usepackage{subfig}
\usepackage{float}
\usepackage{graphicx,tikz}
\usepackage{mathrsfs}
\usepackage[colorlinks=true,citecolor=black,linkcolor=black,urlcolor=black]{hyperref}
\usepackage{hyperref} 
\hypersetup{colorlinks,
    linkcolor=blue, %
    anchorcolor=blue,
    citecolor=blue}
\usepackage{algorithm}

\usepackage{algpseudocode}
\usepackage{chngpage}
\usepackage{mathtools}

\newcommand{\diag}{\mathrm{diag}}

\newtheorem{theorem}{Theorem}[section]
\newtheorem{lemma}[theorem]{Lemma}

\newtheorem{corollary}[theorem]{Corollary}
\newtheorem{definition}[theorem]{Definition}

\newtheorem{proposition}[theorem]{Proposition}

\newtheorem{example}[theorem]{Example}

\theoremstyle{definition}

\DeclareMathOperator{\spec}{sp}

\def\j{\mbox{\boldmath $j$}}

\def\x{\mbox{\boldmath $x$}}
\def\y{\mbox{\boldmath $y$}}
\def\z{\mbox{\boldmath $z$}}
\def\w{\mbox{\boldmath $w$}}

\def\vec0{\mbox{\boldmath $0$}}

\def\B{\mbox{\boldmath $B$}}

\def\I{\mbox{\boldmath $I$}}
\def\J{\mbox{\boldmath $J$}}
\def\L{\mbox{\boldmath $L$}}
\def\M{\mbox{\boldmath $M$}}

\def\O{\mbox{\boldmath $O$}}

\def\SS{\mbox{\boldmath $S$}}

\tikzstyle{vertex}=[circle, draw, inner sep=0pt, minimum size=3pt]
\newcommand{\vertex}{\node[vertex]}

\numberwithin{equation}{section} 
\allowdisplaybreaks  

\def\qed{\hfill$\Box$\vspace{12pt}}



\begin{document}
\title{On the Algebraic Connectivity of Token Graphs\\
and Graphs under Perturbations\thanks{This research has been supported by National Natural Science Foundation of China (No.~12471334, No.~12131013), and Shaanxi Fundamental Science Research Project for Mathematics and Physics (No. 22JSZ009). C. Dalf\'o and M. A. Fiol are funded by AGAUR from the Catalan Government under project 2021SGR00434 and MICINN from the Spanish Government under project PID2020-115442RB-I00.
M. A. Fiol's research is also supported by a grant from the Universitat Polit\`ecnica de Catalunya, reference AGRUPS-2024.}\\
{\small Dedicated to Professor Fuji Zhang on the occasion of his 88th birthday.}}
\author{X. Song$^{a,b}$, C. Dalf\'o$^b$, M. A. Fiol$^c$, S. Zhang$^{a}$\\ 
$^a${\small School of Mathematics and Statistics, Northwestern Polytechnical University}\\
    {\small Xi'an-Budapest Joint Research Center for Combinatorics, Northwestern Polytechnical University} \\
    {\small Xi'an, Shaanxi, P.R. China, {\tt songxd@mail.nwpu.edu.cn, sgzhang@nwpu.edu.cn}}\\
$^b${\small Departament de Matem\`atica, Universitat de Lleida} \\
		{\small Igualada (Barcelona), Catalonia, {\tt{cristina.dalfo@udl.cat}}}\\
$^c${\small Departament de Matem\`atiques, Universitat Polit\`ecnica de Catalunya} \\
    	{\small Barcelona Graduate School, Institut de Matem\`atiques de la UPC-BarcelonaTech (IMTech)}\\
    	{\small Barcelona, Catalonia, {\tt{miguel.angel.fiol@upc.edu}}}
}
\date{}
\maketitle

\begin{abstract}
Given a graph $G=(V,E)$ on $n$ vertices and an integer $k$ between 1 and $n-1$, the $k$-token graph $F_k(G)$ has vertices representing the $k$-subsets of $V$, and two vertices are adjacent if their symmetric difference is the two end-vertices of an edge in $E$. Using the theory of Markov chains of random walks and the interchange process, it was proved that the algebraic connectivities (second smallest Laplacian eigenvalues) of $G$ and $F_k(G)$ coincide, but a combinatorial/algebraic proof has been shown elusive. 
In this paper, we use the latter approach and prove that such equality holds for different new classes of graphs under perturbations, such as extended cycles, extended complete bipartite graphs, kite graphs, and graphs with a cut clique.
Kite graphs are formed by a graph (head) with several paths (tail) rooted at the same vertex and with exciting properties. For instance, we show that the different eigenvalues of a kite graph are also eigenvalues of its perturbed graph obtained by adding edges.
Moreover, as a particular case of one of our theorems, we generalize a recent result of Barik and Verma \cite{bv24} about graphs with a cut vertex of degree $n-1$. 
Along the way, we give conditions under which the perturbed graph $G+uv$, with $uv\in E$, has the same algebraic connectivity as $G$.
\end{abstract}

	

\section{Introduction}

Let $G=(V,E)$ be a graph on $n=|V|$ vertices, with Laplacian matrix $\L=\L(G)$. This is a symmetric positive semidefinite and singular matrix with eigenvalues $\lambda_1(=0)\le\lambda_2\le\cdots\le\lambda_n$. Since its introduction by Fiedler \cite{Fiedler_1973}, the algebraic connectivity of $G$, $\alpha(G)=\lambda_2$, together with its eigenvector $\x$ (or Fiedler vector), has deserved much attention in the literature. See, for instance, the comprehensive survey by Abreu \cite{a07}. Moreover, the behavior of $\alpha(G)$ and $\x$ under graph perturbations has been considered.
Recently, some papers have dealt with the algebraic connectivity of token graphs.
The $k$-token graph $F_k(G)$ has vertices representing the different configurations of $k$ indistinguishable tokens in different vertices of $G$. Two configurations are adjacent if one can be obtained from the other by moving a token along an edge from
its current position to an unoccupied vertex. It was shown that the Laplacian spectrum of a graph $G$ is contained in the spectrum of its $k$-token graph; see 
Dalf\'o, Duque, Fabila-Monroy, Fiol, Huemer, Trujillo-Negrete, and  Zaragoza Mart\'inez \cite{Dalfo_2021}. Moreover, the authors proved that $\alpha(F_k(G)) = \alpha(G)$ for some families of graphs, such as the complete graphs and the paths. From their results, they conjectured that this was always the case.
Later, they realized that this corresponded to the so-called `Aldous spectral graph conjecture' already proved by Caputo, Liggett, and Richthammer \cite{clr10}. The proof was based on Markov chains of random walks and the interchange process. However, until now, a combinatorial/algebraic proof is elusive. Some advances towards such a proof have been recently made; see, for instance, Barik and Verma \cite{bv24}, Dalf\'o and Fiol \cite{Dalfo_2024}, and  Reyes, Dalf\'o, Fiol, and Messegu\'e \cite{Reyes_2025}.  
In this article, our contribution is a further step in this direction by proving the equality $\alpha(F_k(G))=\alpha(G)$ for different new infinite families of graphs.
Besides, along the way, we give the conditions under which the perturbed graph $G+uv$ 
(where $uv$ is a new edge) has the same algebraic connectivity as $G$.

This paper is structured as follows. The following section gives the preliminary results and notations used in this work. The new results are presented in the further sections. Thus, in Section \ref{sec:+edges}, we derive results about the algebraic connectivities of a graph and the same graph after
adding one or more new edges. 
In Section \ref{sec:kite}, we give some results about the so-called kite graphs, that is, with a `head', which is a general graph, and a `tail,' which is a starlike tree (roughly speaking, a tree with paths). 
In the same section, we provide conditions under which the algebraic connectivity does not change when kite graphs are `perturbed' by adding some edges. In Section \ref{sec:cut-clique}, we compute the algebraic connectivity of graphs with a cut clique $K_r$ with maximum degree. This result generalizes a theorem of Barik and Verma \cite{bv24} about the algebraic connectivity and cut vertices. 

\section{Preliminaries}
\label{sec:prelim}
\begin{figure}
    \begin{center}
    \begin{tikzpicture}[scale=1.1,auto,swap]
    \vertex (1) at (0.4,2.25) [fill,label=above:{$1$}]{};
    \vertex (2) at (1.6,2.25) [fill,label=above:{$2$}]{};
    \vertex (3) at (1,1) [fill,label=left:{$3$}]{};
    \vertex (4) at (1,-0.2) [fill,label=left:{$4$}]{};
    \vertex (5) at (1,-1.4) [fill,label=left:{$5$}]{};
    \draw[line width=0.6pt](1)--(3);
    \draw[line width=0.6pt](2)--(3);
    \draw[line width=0.6pt](3)--(4);
    \draw[line width=0.6pt](4)--(5);
    \vertex (12) at (6.5,2.55) [fill,label=above:{$12$}]{};
    \vertex (34) at (6.5,1.3) [fill,label=above:{$34$}]{};
    \vertex (35) at (6.5,0.1) [fill]{};
    \node () at (6.6,-0.05) [label=above left:{$35$}]{};
    \vertex (45) at (6.5,-1.5) [fill,label=left:{$45$}]{};
    \vertex (13) at (5.2,2.05) [fill,label=left:{$13$}]{};
    \vertex (14) at (5.2,0.8) [fill,label=left:{$14$}]{};
    \vertex (15) at (5.2,-0.4) [fill,label=left:{$15$}]{};
    \vertex (23) at (7.8,2.05) [fill,label=right:{$23$}]{};
    \vertex (24) at (7.8,0.8) [fill,label=right:{$24$}]{};
    \vertex (25) at (7.8,-0.4) [fill,label=right:{$25$}]{};
    \draw[line width=0.6pt](12)--(13);
    \draw[line width=0.6pt](12)--(23);
    \draw[line width=0.6pt](23)--(24);
    \draw[line width=0.6pt](13)--(14);
    \draw[line width=0.6pt](14)--(34);
    \draw[line width=0.6pt](34)--(24);
    \draw[line width=0.6pt](24)--(25);
    \draw[line width=0.6pt](34)--(35);
    \draw[line width=0.6pt](14)--(15);
    \draw[line width=0.6pt](35)--(15);
    \draw[line width=0.6pt](35)--(25);
    \draw[line width=0.6pt](35)--(45);
    \end{tikzpicture}
    \caption{The graph $Y$ and its $2$-token graph $F_2(Y)$.}\label{fig:Y+F2(Y)}
    \end{center}
\end{figure}
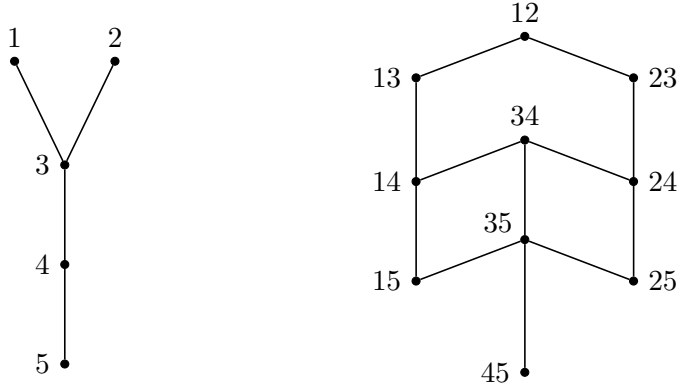

Let $G=(V,E)$ be a graph on $n$ vertices. Let $V'\subset V$. Denote by $G[V']$ the {\em induced subgraph} of $G$ whose vertex set is $V'$ and whose edge set consists of all edges of $G$ with both ends in $V'$. If $G'=(V',E')$ is an induced subgraph of $G$, the {\em degree} of $G'$, denoted by $d_{G}(G')$, is the number of edges between $V'$ and $V\setminus V'$. In particular, the degree of a vertex $v$ is denoted by $d_G(v)$. 
We denote by $G_1\cup G_2$ the {\em union} of two simple graphs $G_1=(V_1,E_1)$ and $G_2=(V_2,E_2)$, with vertex set $V_1\cup V_2$ and edge set $E_1\cup E_2$. If $G_1$ and $G_2$ are disjoint, we refer to their union as a {\em disjoint union}.
Given an integer $k$ such that $1\le k\le n-1$, the {\em $k$-token graph} of $G$, denoted by $F_k(G)$, has ${n\choose k}$ vertices representing the configurations of $k$ indistinguishable tokens
placed at distinct vertices of $G$. Moreover,  two configurations are adjacent whenever one
configuration can be reached from the other by moving one token along an edge of $G$ from its current position to an unoccupied vertex, see 
Fabila-Monroy, Flores-Peñaloza, Huemer, Hurtado, Urrutia, and Wood \cite{ffhhuw12}.
As an example, Figure \ref{fig:Y+F2(Y)} shows the graph $Y$ on 5 vertices, and its 2-token graph $F_2(Y)$ with ${5\choose 2}=10$ vertices.


For a matrix $\M$ with order $n$, let $\lambda_1(\M)\le \lambda_2(\M)\le \cdots \le \lambda_n(\M)$ be its eigenvalues. 
Denote by $\Phi(\M)$ the characteristic polynomial of $\M$.
Let $\L=\L(G)$ be the Laplacian matrix of $G$. This matrix is positive semidefinite and singular, with eigenvalues $(0=)\lambda_1 \le \lambda_2\le \cdots\le \lambda_n$. Its second smallest eigenvalue $\lambda_2$ is
known as the {\em algebraic connectivity} of $G$ (see Fiedler \cite{Fiedler_1973}), and we denote it by $\alpha(G)$.
An eigenvector corresponding to the algebraic connectivity is called a {\em Fiedler vector}. 

Let $\j$ be the all-1 (column) vector. Given a graph $G$ of order $n$, we define a vector $\x$ to be an \emph{embedding} of $G$ if $\x \in W_n$, where $W_n:=\{\y:\y^{\top}\j = 0\}$. 
For a vertex $v\in V$, the entry of $\x$ corresponding to $v$ is denoted by $\x_v$.
The value $\frac{\x^{\top}\L(G)\x}{\x^{\top}\x}$ is known as the \emph{Rayleigh quotient}, and 
\begin{equation}
 \alpha(G)\le \frac{\x^{\top}\L(G)\x}{\x^{\top}\x}=\frac{\sum_{uv\in E(G)}(\x_u-\x_v)^2}{\x^{\top}\x}\label{eq-rayleigh}   
\end{equation}
for any vector $\x\in W_n$ and $\x\ne \vec0$, where equality holds if and only if $\x$ is a Fiedler vector of $G$.

Given some integers $n$ and $k$ (with $k\in [n]=\{1,2,\ldots\}$),  the $(n,k)$-\emph{binomial matrix} $\B$  is a ${n \choose k}\times n$ matrix whose rows are the characteristic vectors of the $k$-subsets of $[n]$ in a given order. Thus, if the $i$-th $k$-subset is $A$, then
$$
(\B)_{ij}=
\left\lbrace
\begin{array}{ll}
	1 & \mbox{if } j\in A,\\
	0 & \mbox{otherwise.}
\end{array}
\right.
$$

 

Denote by $\O$ the all-$0$ matrix with the appropriate dimension.
Let $\I_n$ and $\J_n=\j\j^{\top}$ be the identity matrix and the all-$1$ matrix with order $n$, respectively.

Now, let us describe a pair of results on the interlacing of graph eigenvalues.

\begin{lemma}[Horn and Johnson \cite{Horn_2013}]
\label{le-principal_submatrix}
Let $\M'$ be a Hermitian square matrix of order $n-1$, $\z\in \mathbb{C}^{n-1}$, and $a\in \mathbb{R}$, and let 
$$
\M=\left({\begin{array}{cccc}
        \M'&\z\\
        \z^*&a
        \end{array}}
        \right),
$$
where $\z^*$ is the conjugate transpose of $\z$. 
Then, the eigenvalues of $\M'$ interlace the eigenvalues of $\M$. That is, we have the inequalities
\[
\lambda_1(\M)\le \lambda_1(\M')\le \lambda_2(\M)\le \cdots\le \lambda_{n-1}(\M)\le \lambda_{n-1}(\M')\le \lambda_n(\M),
\]
in which $\lambda_i(\M')=\lambda_{i+1}(\M)$ if and only if there is a nonzero vector $\w \in \mathbb{C}^{n-1}$ such that $\M'\w=\lambda_i(\M')\w$, with $\z^*\w=0$, and $\M'\w=\lambda_{i+1}(\M)\w$. 
\end{lemma}
\begin{lemma}[Horn and Johnson \cite{Horn_2013}]
\label{le-matrix_add}
    Let $\M_1$ and $\M_2$ be two Hermitian matrices of order $n$. Then,
    \[
    \lambda_i(\M_1)+\lambda_1(\M_2)\le \lambda_i(\M_1+\M_2)\le \lambda_i(\M_1)+\lambda_n(\M_2)
    \]
    for $i=1,\ldots,n$, with equality in the upper bound if and only if there is nonzero vector $\x$ such that $\M_1\x=\lambda_i(\M_1)\x$, $\M_2\x=\lambda_n(\M_2)\x$, and $(\M_1+\M_2)\x=\lambda_i(\M_1+\M_2)\x$; equality in the lower bound occurs if and only if there is a nonzero vector $\x$ such that $\M_1\x=\lambda_i(\M_1)\x$, $\M_2\x=\lambda_1(\M_2)\x$, and $(\M_1+\M_2)\x=\lambda_i(\M_1+\M_2)\x$.
\end{lemma}

Let us now show some results concerning a graph after adding a new edge $uv$ to it.

\begin{lemma}[Cvetkovi\'c, Doob, and Sachs \cite{Cvetkovic_1980}]
\label{le-interlacing}
Let $G$ be a graph and $G'=G+uv$. Then, the Laplacian eigenvalues of $G$ and $G'$ interlace, that is, 
\[
0=\lambda_1(G)= \lambda_1(G')\le \lambda_2(G)\le \lambda_2(G')\le \cdots\le \lambda_n(G) \le \lambda_{n}(G').
\]
\end{lemma}

\begin{lemma}[Xue, Lin, and Shu \cite{Xue_2019}]
\label{le-Xue}
Let $G$ be a graph and $G'=G+uv$. If $\alpha(G')=\alpha(G)$, then there exists a Fiedler vector $\x$ of $G$ such that $\x_u=\x_v$. 
\end{lemma}
\begin{lemma}[Merris \cite{Merris_1998}]
\label{le-Merris}
Let $G$ be a graph and $G'=G + uv$ or $G'=G - uv$. Let $\lambda$ be a Laplacian eigenvalue of $G$ corresponding to an eigenvector $\x$. If $\x_u=\x_v$, then $\lambda$ is also a Laplacian eigenvalue of $G'$ with eigenvector  $\x$, where $G'$ is the graph obtained from $G$ by deleting or adding an edge $uv$ depending on whether or not it is an edge of $G$. 
\end{lemma}

From Lemma \ref{le-Merris}, we 
get that $\lambda=\alpha(G)$ is also a Laplacian eigenvalue of $G+uv$ corresponding to $\x$ if $\x_u=\x_v$. Then, it must be $\alpha(G+uv)\le \alpha(G)$ and, since $\alpha(G+uv)\ge \alpha(G)$, we  conclude that
$\alpha(G)=\alpha(G+uv)$.
Combining Lemmas \ref{le-Xue} and \ref{le-Merris}, we can obtain the following lemma, which gives a \textbf{necessary and sufficient condition} for $\alpha(G)=\alpha(G+uv)$.

\begin{lemma}
\label{le-adding_edge_iff}
Let $G$ be a graph on $n$ vertices.
Then, $\alpha(G)=\alpha(G+uv)$ if and only if there exists a Fiedler vector $\x$ of $G$ such that $\x_u=\x_v$. 
\end{lemma}
Of course, by applying repeatedly Lemma \ref{le-adding_edge_iff},  if $G'$
is obtained from $G$ by adding $r$ edges $u_iv_i$, we have that  $\alpha(G)=\alpha(G')$ if and only if there exists a Fiedler vector $\x$
of $G$ such that $\x_{u_i}=\x_{v_i}$ for $i=1,\ldots,r$.

Some preliminary results related to token graphs are presented below.
\begin{lemma}[Dalf\'o, Duque, Fabila-Monroy, Fiol, Huemer, Trujillo-Negrete, and
Zaragoza Mart\'{\i}ınez \cite{Dalfo_2021}]
\label{le-Dalfo2021}
Consider a graph $G(\cong F_1(G))$ and its $k$-token graph $F_k(G)$, with corresponding Laplacian matrices $\L(G)$ and $\L(F_k(G))$, and the $(n,k)$-binomial matrix $\B$. The following statements hold:
\begin{itemize}
\item[$(i)$] If $\x$ is a $\lambda$-eigenvector of $\L(G)$, then $\B\x$ is a $\lambda$-eigenvector of $\L(F_k(G))$.
\item[$(ii)$] If $\w$ is a $\lambda$-eigenvector of $\L(F_k(G))$ and $\B^{\top}\w\ne 0$, then $\B^{\top}\w$ is a $\lambda$-eigenvector of $\L(G)$.
\end{itemize}
\end{lemma}

As a consequence of Lemma \ref{le-Dalfo2021}, we have the following result.

\begin{theorem}[Dalf\'o, Duque, Fabila-Monroy, Fiol, Huemer, Trujillo-Negrete, and
Zaragoza Mart\'{\i}ınez \cite{Dalfo_2021}]
\label{th-Dalfo2021}
    Let $G$ be a graph on $n$ vertices and let $k$ be an integer such that $1\le k\le n-1$. Then, the spectrum of $G$ is contained in the spectrum of its $k$-token graph $F_k(G)$. 
\end{theorem}

\begin{theorem}[Dalf\'o and Fiol \cite{Dalfo_2024}, Dalf\'o, Duque, Fabila-Monroy, Fiol, Huemer, Trujillo-Negrete, and
Zaragoza Mart\'{\i}ınez \cite{Dalfo_2021}, Reyes, Dalf\'o, Fiol, and Messegu\'e \cite{Reyes_2025}]
\label{th-Dalfo2024}
For each of the following classes of graphs, the algebraic connectivity of a token graph $F_k(G)$ satisfies the following statements.
\begin{itemize}
\item[$(i)$]
Let $T_n$ be a tree on $n$ vertices. Then,
$\alpha(F_k(T_n))=\alpha(T_n)$ for every $n$ and $k=1,\ldots,n-1$.
\item[$(ii)$]
Let $G$ be a graph such that $\alpha(F_k(G))=\alpha(G)$. Let $T_G$ be a graph in which each vertex of $G$ is the root vertex of some (possibly empty) tree. Then, 
$\alpha(F_k(T_G))=\alpha(T_G)$.
\item[$(iii)$]
Let $G=K_{n_1,n_2}$ be a complete bipartite graph on $n=n_1+n_2$ vertices, with $n_1\le n_2$. Then, $\alpha(F_k(G))=\alpha(G)=n_1$ for every $n_1,n_2$ and $k=1,\ldots,n-1$.
\item[$(iv)$]
Let $G=K_n$ be a complete graph on $n$ vertices. Then, $\alpha(F_k(G))=\alpha(G)=n$ for every $k=1,\ldots,n-1$.
\item[$(v)$]
Let $C_n$ be a cycle on $n$ vertices. Then, $\alpha(F_k(C_n))=\alpha(C_n)$ for $k=1,2$.
\end{itemize}
\end{theorem}
\begin{lemma}[Barik and Verma \cite{bv24}]
Let $G$ be a graph on $n\ge 4$ vertices, and $H$ be the graph formed by adding a
pendant vertex (say $v$, with the corresponding edge) to graph $G$. Then, for any integer k such that $2\le k\le n/2$,
$$
\alpha(F_k(H)) \le  \min\{\alpha(F_{k-1}(G)) + 1, \alpha(F_k(G)) + 1\}.
$$
\end{lemma}
\begin{theorem}[Barik and Verma \cite{bv24}]
\label{theBarik}
Let $G$ be a graph on $n$ vertices and $v$ be a cut vertex in $G$ such that $d_G(v)=n-1$. Then, for any integer $k$ such that $2\le k\le \frac{n}{2}$,
$$
\alpha(F_k(G))=\alpha(G)=1.
$$
\end{theorem}

\section{Adding edges}
\label{sec:+edges}
In this section and the following ones, we present the new results of the paper. We begin by deriving results about the algebraic connectivities of a graph and the same graph after adding a new edge $uv$.


\subsection{A basic result}

\begin{theorem}
\label{th-add_edge}
Let $G=(V,E)$ be a graph with order $n$ such that, for some $k$ with $2\le k\le n/2$, we have $\alpha(F_k(G))=\alpha(G)$. Consider adding 
an edge
$uv$, for $u,v\in V$, 
getting the new graph $G+uv$. 
If $\alpha(G+uv)=\alpha(G)$, then $\alpha(F_k(G+uv))=\alpha(G+uv)$.
\end{theorem}
\begin{proof}
    Note that $F_k(G)$ is a spanning subgraph of $F_k(G+uv)$ with  
    \[
    E(F_k(G+uv))\backslash E(F_k(G))=\{A_{r}A_{s} : A_{r}=\{u,u_1,\ldots,u_{k-1}\}, A_{s}=\{v,u_1,\ldots,u_{k-1}\}\},
     \]
     where $u_1,\ldots,u_{k-1}\in V(G)\backslash \{u,v\}$.  
    Then, by Lemma \ref{le-interlacing},
    \begin{equation}
        \alpha(F_k(G+uv)){\ge}\alpha(F_k(G)).\label{eqthe2'}
    \end{equation}
    Since $\alpha(G+uv)=\alpha(G)$, by Lemma \ref{le-adding_edge_iff}, there exists a Fiedler vector $\x$ of $G$ such that $\x_{u}=\x_{v}$. 
    Let $\y=\B\x$, where $\B$ is the $(n,k)$-binomial matrix. 
    It follows from the hypothesis $\alpha(F_k(G))=\alpha(G)$ and Lemma \ref{le-Dalfo2021}($i$) that $\y$ is a Fiedler vector of $F_k(G)$. Moreover, observe that
    $$
    \y_{A_r}-\y_{A_s}=\x_u+\sum_{i=1}^{k-1}\x_{u_i}-\left(\x_v+\sum_{i=1}^{k-1}\x_{u_i}\right)=\x_u-\x_v=0.
    $$
    Hence, we get 
    \begin{align}
         \alpha(F_k(G))=\frac{\y^{\top}\L(F_k(G))\y}{\y^{\top}\y}&=\frac{\sum_{A_rA_s\in E(F_k(G))}(\y_{A_r}-\y_{A_s})^2}{\y^{\top}\y}\notag\\
       &=\frac{\sum_{A_rA_s\in E(F_k(G+uv))}(\y_{A_r}-\y_{A_s})^2}{\y^{\top}\y}\notag\\
       &\ge \alpha(F_k(G+uv)).\label{eqthe3'}
    \end{align}
    Thus, from \eqref{eqthe2'}, 
    \eqref{eqthe3'}, and the hypothesis, we conclude that $\alpha(F_k(G+uv))=\alpha(F_k(G))=\alpha(G)=\alpha(G+uv)$.
    
    Alternatively, we have
    \begin{equation*}
        \alpha(F_k(G+uv))\ge\alpha(F_k(G))=\alpha(G)=\alpha(G+uv).
        \label{eq1'}
    \end{equation*}
    Then, the result follows since, by Theorem \ref{th-Dalfo2021}, it must be  $\alpha(F_k(G+uv))\le \alpha(G+uv)$.\qed
\end{proof}

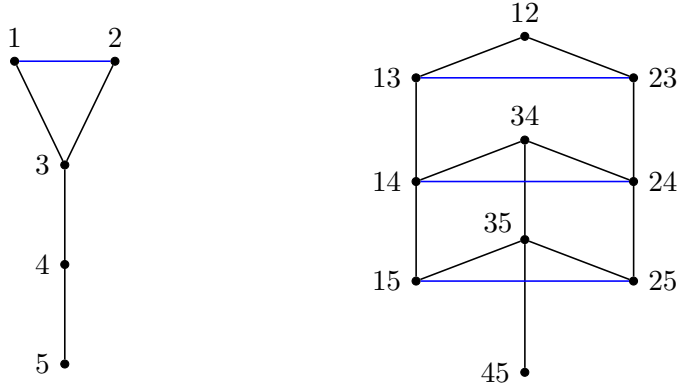
\begin{figure}
    \begin{center}
    \begin{tikzpicture}[scale=1.1,auto,swap]
    \vertex (1) at (0.4,2.25) [fill,label=above:{$1$}]{};
    \vertex (2) at (1.6,2.25) [fill,label=above:{$2$}]{};
    \vertex (3) at (1,1) [fill,label=left:{$3$}]{};
    \vertex (4) at (1,-0.2) [fill,label=left:{$4$}]{};
    \vertex (5) at (1,-1.4) [fill,label=left:{$5$}]{};
    \draw[line width=0.6pt](1)--(3);
    \draw[blue,line width=0.6pt](1)--(2);
    \draw[line width=0.6pt](2)--(3);
    \draw[line width=0.6pt](3)--(4);
    \draw[line width=0.6pt](4)--(5);
    \vertex (12) at (6.5,2.55) [fill,label=above:{$12$}]{};
    \vertex (34) at (6.5,1.3) [fill,label=above:{$34$}]{};
    \vertex (35) at (6.5,0.1) [fill]{};
    \node () at (6.6,-0.05) [label=above left:{$35$}]{};
    \vertex (45) at (6.5,-1.5) [fill,label=left:{$45$}]{};
    \vertex (13) at (5.2,2.05) [fill,label=left:{$13$}]{};
    \vertex (14) at (5.2,0.8) [fill,label=left:{$14$}]{};
    \vertex (15) at (5.2,-0.4) [fill,label=left:{$15$}]{};
    \vertex (23) at (7.8,2.05) [fill,label=right:{$23$}]{};
    \vertex (24) at (7.8,0.8) [fill,label=right:{$24$}]{};
    \vertex (25) at (7.8,-0.4) [fill,label=right:{$25$}]{};
    \draw[line width=0.6pt](12)--(13);
    \draw[line width=0.6pt](12)--(23);
    \draw[line width=0.6pt](23)--(24);
    \draw[line width=0.6pt](13)--(14);
    \draw[line width=0.6pt](14)--(34);
    \draw[line width=0.6pt](34)--(24);
    \draw[line width=0.6pt](24)--(25);
    \draw[line width=0.6pt](34)--(35);
    \draw[line width=0.6pt](14)--(15);
    \draw[line width=0.6pt](35)--(15);
    \draw[line width=0.6pt](35)--(25);
    \draw[line width=0.6pt](35)--(45);
    \draw[blue,line width=0.6pt](13)--(23);
    \draw[blue,line width=0.6pt](14)--(24);
    \draw[blue,line width=0.6pt](15)--(25);
    \end{tikzpicture}
    \caption{The graph $Y+12$ and its 2-token graph $F_2(Y+12)$.}\label{fig:F2(Y+12)}
    \end{center}
\end{figure}

\begin{example}
The graph $Y=(V,E)$ of Figure \ref{fig:Y+F2(Y)} has algebraic connectivity $\alpha(Y)=0.5188$ with corresponding eigenvector $\x=(-0.5969,-0.5969,-2.8772,0.4812,1)^{\top}$. Then, since $\x_1=\x_2$, its `extended' graph $Y+12$ has the same algebraic connectivity
$\alpha(Y+12)=\alpha(Y)$, with the same eigenvector $\x$.
This graph is shown in Figure \ref{fig:F2(Y+12)} together with its 2-token graph $F_2(Y+12)$. Notice that $F_2(Y)$ is a spanning subgraph of $F_2(Y+12)$, where the `new' edges induced by $12$ are $A_rA_s\in \{13\sim 23, 14\sim 24, 15\sim 25\}$. Since $\alpha(F_2(Y))=\alpha(Y)$, Theorem \ref{th-add_edge} implies that  $\alpha(F_2(Y+12))=\alpha(Y+12)=\alpha(Y)=0.5188$. 
\end{example}
%

\subsection{Extended graphs with pendant vertices}

From the result in Theorem \ref{th-add_edge},
it is natural to consider graphs satisfying $\alpha(G)=\alpha(G+uv)$ for some edge $uv$.
A family of such graphs is given in the following Lemma \ref{le-Shao2008}, whose statement
can be made more explicit by first computing the value of a particular eigenvalue. With this aim, 
given a vertex subset $V'\subseteq V$, let $\L_{[V']}(G)$ denote the principal submatrix of $\L(G)$ whose rows and columns correspond to $V'$. 
When $G$ is a path $P_{r+1}$ with vertices $u_0,u_1,\ldots,u_r$ and $V_r=\{u_1,\ldots,u_r\}$, $\L_{[V_r]}(G)$ is the $r\times r$ tridiagonal matrix
$$
  \L_{[V_r]}(P_{r+1}) =
  {\scriptsize\left(
  \begin{array}{cccccc}
  2 & -1 & 0 & 0 & 0 & 0\\
  -1 & 2 & -1 & 0 & 0 & 0\\
  0 & -1 & \ddots & \ddots & 0 & 0\\
  0 & 0 & \ddots & \ddots & -1 & 0\\
  0 & 0 & 0 & -1 & 2 & -1 \\
  0 & 0 & 0 & 0 & -1 & 1 
  \end{array}
  \right)},
  $$
with eigenvalues
\begin{equation}
\label{eq:theta-k}
\theta_k=2+2\cos\left(\frac{2k\pi}{2r+1}\right),\quad\mbox{for $k=1,\ldots,r$}.
\end{equation}
(See Yueh \cite[Th. 1]{y05}). 
The minimum eigenvalue is obtained when $k=r$ (since $\theta_k$ is decreasing with $k$), so we have that
 $\lambda_1(\L_{[V_r]}(P_{r+1}))=\theta_r$, and the following result of 
 Shao, Guo, and Shan \cite{Shao_2008} reads as follows.
 
\begin{lemma}[Shao, Guo, and Shan \cite{Shao_2008}]
\label{le-Shao2008}
Let $v$ be a vertex in a connected graph $H$ and suppose that $s(\ge 2)$ new paths (with equal length $r$) $vv_{i1}v_{i2}\cdots v_{ir}$ (for $i=1,\ldots,s$ and $r\ge 1$) are attached to $H$ at $v$ to form a new graph $G$. 
Let $G^+$ be the graph obtained from $G$ by arbitrarily adding edges among the vertices $v_{1j},v_{2j},\ldots,v_{sj}$ for any given $j=1,\ldots,r$. 
If $\alpha(G)\neq 2+2\cos\left(\frac{2r\pi}{2r+1}\right)$,
then $\alpha(G^+)=\alpha(G)$.
\end{lemma}

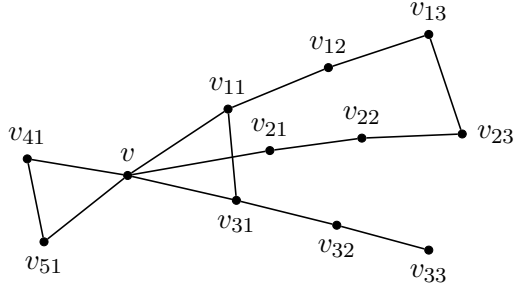
\begin{figure}
    \begin{center}
    \begin{tikzpicture}[scale=1.1,auto,swap]
    \vertex (0) at (1,1) [fill,label=above:{$v$}]{};
    \vertex (11) at (2.2,1.8) [fill,label=above:{$v_{11}$}]{};
    \vertex (12) at (3.4,2.3) [fill,label=above:{$v_{12}$}]{};
    \vertex (13) at (4.6,2.7) [fill,label=above:{$v_{13}$}]{};
    \vertex (21) at (2.7,1.3) [fill,label=above:{$v_{21}$}]{};
    \vertex (22) at (3.8,1.45) [fill,label=above:{$v_{22}$}]{};
    \vertex (23) at (5,1.5) [fill,label=right:{$v_{23}$}]{};
    \vertex (31) at (2.3,0.7) [fill,label=below:{$v_{31}$}]{};
    \vertex (32) at (3.5,0.4) [fill,label=below:{$v_{32}$}]{};
    \vertex (33) at (4.6,0.1) [fill,label=below:{$v_{33}$}]{};
    \vertex (41) at (-0.2,1.2) [fill,label=above:{$v_{41}$}]{};
    \vertex (51) at (0,0.2) [fill,label=below:{$v_{51}$}]{};
    \draw[line width=0.6pt](0)--(11);
    \draw[line width=0.6pt](0)--(21);
    \draw[line width=0.6pt](0)--(31);
    \draw[line width=0.6pt](0)--(41);
    \draw[line width=0.6pt](0)--(51);
    \draw[line width=0.6pt](11)--(12);
    \draw[line width=0.6pt](12)--(13);
    \draw[line width=0.6pt](13)--(23);
    \draw[line width=0.6pt](11)--(31);
    \draw[line width=0.6pt](21)--(22);
    \draw[line width=0.6pt](22)--(23);
    \draw[line width=0.6pt](31)--(32);
    \draw[line width=0.6pt](32)--(33);
    \draw[line width=0.6pt](41)--(51);
    \end{tikzpicture}
    \caption{An example of $G^+$ where $H=K_3$ and $s=r=3$ in Lemma \ref{le-Shao2008}.}\label{fig:example}
    \end{center}
\end{figure}

\begin{example}
A simple example of Lemma \ref{le-Shao2008} is obtained when $G=Y$, the graph of Figure \ref{fig:Y+F2(Y)}, where $s=2$ and $r=1$. Then, as $\alpha(Y)=0.5188\neq \theta_1=1$, the graph $G^+=Y+12$ of Figure \ref{fig:F2(Y+12)} satisfies $\alpha(G^+)=\alpha(G)$.
Moreover, Figure \ref{fig:example} is an example when the statement of Lemma \ref{le-Shao2008} does not hold because the necessary conditions fail. In the latter graph, we have $\alpha(G^+)=0.2679$, which is different from $\alpha(G^+-v_{13}v_{23})=0.1981$. This is because $\alpha(G)=\alpha(G^+-\{v_{13}v_{23}, v_{11}v_{31}\})=\lambda_1(\L_{[V_r]}(P_{r+1}))=\theta_r=0.1981$ (see Table \ref{tab1} for $r=3$).
\begin{table}
\begin{center}
{\small
\begin{tabular}{|c||cccccccccc|}
\hline
$r$ & 1 & 2 & 3 & 4 & 5 & 6 & 7 & 8 & 9 & 10\\
\hline
$\theta_r$ & 1  & 0.3820  & 0.1981 & 0.1206  & 0.0810  & 0.0581 & 0.0437  & 0.0341 & 0.0273  & 0.0223\\
\hline
\end{tabular}
}
\end{center}
\caption{The values of $\theta_r=\lambda_1(\L_{[V_r]}(P_{r+1}))$ for different values of $r$.}
\label{tab1}
\end{table}
\end{example}

It has been shown in Kirland \cite{Kirkland_2000} that, for a connected graph $G$ with $n$ vertices and a cut vertex $v$, we have $\alpha(G)\le 1$ and the equality holds if and only if $d_G(v)=n-1$. We get the following theorem by combining this with 
previous results.
\begin{theorem}
\label{th-addedge-pendant}
    Let $G$ be a graph on $n$ vertices. Suppose that $v_1,v_2,\ldots,v_s$ (for $s\ge2$) are $s$ pendant vertices of $G$ adjacent to a common vertex $v$ with $d_G(v)\ne n-1$. 
    Let $G^+$ be a graph obtain from $G$ by adding any $t$ $\left(\mbox{for } 0\le t\le \frac{s(s-1)}{2}\right)$ edges among $v_1,v_2,\ldots,v_s$. 
    If $\alpha(F_k(G))=\alpha(G)$, then $\alpha(F_k(G^+))=\alpha(G^+)$. 
\end{theorem}
\begin{proof}
    Note that $v$ is a cut vertex of $G$. Since $d_G(v)\ne n-1$, we have $\alpha(G)<1=\theta_1$. It follows from Lemma \ref{le-Shao2008} and Theorem \ref{th-add_edge} that the results hold.\qed
\end{proof}



The following corollary holds directly from Theorems \ref{th-Dalfo2024}$(i)$ and \ref{th-addedge-pendant}.
\begin{corollary}
    Let $T$ be a tree on $n$ vertices. Suppose that $v_1,v_2,\ldots,v_s$ (for $s\ge2$)  are $s$ pendant vertices of $T$ adjacent to a common vertex $v$ with $d_T(v)\ne n-1$. 
    Let $T^+$ be a graph obtain from $T$ by adding any $t$ $\left(\mbox{for } 0\le t\le \frac{s(s-1)}{2}\right)$ edges among $v_1,v_2,\ldots,v_s$. 
    Then, $\alpha(F_k(T^+))=\alpha(T^+)$.
\end{corollary}

\subsection{Extended cycles}
Let $G=(V,E)=C_n$ be a cycle with $n$ vertices, $v_0,v_1,\ldots v_{n-1}$, and let $C^+_n=(V,E^+)$ be a graph obtained from $C_n$ by adding some edges in  $E_0=\{v_iv_j:i+j=\nu\}$, where $\nu=n$ if $n$ is odd and  $\nu\in \{n,n-1\}$ otherwise, so that $E^+=E\cup E_0$.
Guo, Zhang, and Yu \cite{Guo_2018} showed 
 that there exists a Fiedler vector $\x$ of $C_n$ such that $\x_{v_i}=\x_{v_{n-i}}$ (for $i=1,\ldots,(n-1)/2$), and an integer $s$ (with $1\le s\le (n-3)/2$) such that $\x_{v_0}>\x_{v_1}>\cdots>\x_{v_s}>0\ge \x_{v_{s+1}}>\cdots>\x_{v_{\frac{n-1}{2}}}$ for odd $n$. (For even $n$, there is a similar result). Consequently, they proved that 
adding the edges in $E_0$ does not change the algebraic connectivity: $\alpha(C_n)=\alpha(C^+_n)$.
Combining this fact with Theorems \ref{th-add_edge} and \ref{th-Dalfo2024}$(v)$, we obtain the following result. 
\begin{theorem}
    Let $C_n$ be a cycle with $n$ vertices, and $C^+_n$ a graph obtained from $C_n$ by adding some edges in $E_0=\{v_iv_j:i+j=\nu\}$. 
    Then, $\alpha(F_{k}(C^+_n))=\alpha(C^+_n)$ for $k=1,2$.
\end{theorem}

\subsection{Extended complete bipartite graphs}
Now, we consider adding new edges to complete bipartite graphs.
\begin{theorem}
\label{th-bip-addedge1}
    Let $K_{n_1,n_2}$ be a complete bipartite graph on $n=n_1+n_2$ vertices, with bipartition $(X,Y)$ such that $|X|=n_1\le n_2=|Y|$. Let $K^+_{n_1,n_2}$ be a graph obtained from $K_{n_1,n_2}$ by adding any $t$ $\left(\mbox{for } 0\le t\le \frac{n_1(n_1-1)}{2}\right)$ edges among $X$. Then, $\alpha(F_k(K^+_{n_1,n_2}))=\alpha(K^+_{n_1,n_2})=n_1$ for every $n_1,n_2$ and $k=1,\ldots,n-1$.
\end{theorem}
\begin{proof}
    Denote by $K^*_{n_1,n_2}$ the graph obtained from $K_{n_1,n_2}$ by adding all edges among $X$, that is, $G[X]\cong K_{n_1}$.
    It can be shown that the spectrum of $K^*_{n_1,n_2}$ is $\{0,n_1^{n_2-1},n^{n_1}\}$.
    Thus, $\alpha(K^*_{n_1,n_2})=n_1$. 
    It follows from Lemma \ref{le-interlacing} that 
    \begin{equation*}
        \alpha(K^*_{n_1,n_2})\ge \alpha(K^+_{n_1,n_2})\ge\alpha(K_{n_1,n_2}).
    \end{equation*}
    Note that $\alpha(K^*_{n_1,n_2})=\alpha(K_{n_1,n_2})=n_1$. We obtain $\alpha(K^+_{n_1,n_2})=n_1=\alpha(K_{n_1,n_2})$.
    Then, the result follows from Theorems \ref{th-Dalfo2024}$(iii)$ and \ref{th-add_edge}.\qed
\end{proof}
How about adding edges among the vertices in $Y$ in the complete bipartite graph? We consider adding all edges among $Y$ based on the following lemma. 
For a vertex $v\in V(G)$, let $S_v:=\{A\in V(F_k(G)):v\in A\}$ and $S_v':=\{B\in V(F_k(G)):v\notin B\}$. 
Let $H_v$ and $H_v'$ be the subgraphs of $F_k(G)$ induced by $S_v$ and $S_v'$, respectively. Note that $H_v\cong F_{k-1}(G-v)$ and $H_v'\cong F_k(G-v)$. 
\begin{lemma}[Dalf\'o and Fiol \cite{Dalfo_2024}]
    \label{le-embedding}
    Given a vertex $v\in V(G)$ and an eigenvector $\w$ of $F_k(G)$ such that $\B^{\top}\w=\vec0$, let $\y_v:=\w|_{S_v}$ and $\y_v':=\w|_{S_v'}$.
    Then, $\y_v$ and $\y_v'$ are embeddings of $H_v$ and $H_v'$, respectively.
\end{lemma}
\begin{theorem}
\label{th-bip-addedge2}
    Let $G=K_{n_1,n_2}$ be a complete bipartite graph on $n=n_1+n_2$ vertices, with bipartition $(X,Y)$ such that $2\le n_1\le n_2$, where $|X|=n_1$ and $|Y|=n_2$. Let $K^*_{n_1,n_2}$ be a graph obtained from $K_{n_1,n_2}$ by adding all edges among $Y$. Then, $\alpha(F_k(K^*_{n_1,n_2}))=\alpha(K^*_{n_1,n_2})=n_2$ for every $n_1,n_2$ and $k=1,\ldots,n-1$.
\end{theorem}

\begin{proof}
    Note that the spectrum of $K^*_{n_1,n_2}$ is $\{0,n_2^{n_1-1},n^{n_2}\}$.
    Thus, $\alpha(K^*_{n_1,n_2})=n_2$.
    Recall that $\alpha(F_k(K^*_{n_1,n_2}))\le \alpha(K^*_{n_1,n_2})=n_2$. 
    Let $\w$ be a Fiedler vector of $F_k(K^*_{n_1,n_2})$. Without loss of generality, suppose that $\w^{\top}\w=1$. 
    It follows from Lemma \ref{le-Dalfo2021}($ii$) that, if $\B^{\top}\w\ne\vec0$, then $\alpha(F_k(K^*_{n_1,n_2}))$ is also an eigenvalue of $K^*_{n_1,n_2}$, so that $\alpha(F_k(K^*_{n_1,n_2}))=\w^{\top}\L(F_k(K^*_{n_1,n_2}))\w\ge \alpha(K^*_{n_1,n_2})=n_2$.
    Then, it suffices to show $\alpha(F_k(K^*_{n_1,n_2}))\ge n_2$ when $\B^{\top}\w=\vec0$.
    We prove it by induction on $n_1(\ge2)$ and $k$. For $k=1$, the claim holds, since $F_1(K^*_{n_1,n_2})\cong K^*_{n_1,n_2}$. 
    For $n_1=2$, note that the claim is true if $k=1$. Then, we assume that $\alpha(F_{k-1}(K^*_{n_1,n_2}))\ge n_2.$
    Let $v$ be a vertex in $X$. 
    Thus, $G-v\cong K_{n_2+1}$ is a complete graph.
    Let $\y_v:=\w|_{S_v}$ and $\y_v':=\w|_{S_v'}$ such that, by Lemma \ref{le-embedding}, are embeddings of $H_v\cong F_{k-1}(G-v)$ and $H_v'\cong F_k(G-v)$. Then,
from (\ref{eq-rayleigh}), 
the hypothesis, and Lemma \ref{th-Dalfo2024}$(iv)$, we have 
\[
\frac{\sum_{(A,B)\in E(H_v)}\left((\y_v)_A-(\y_v)_B\right)^2}{\sum_{A\in V(H_v)}((\y_v)_A)^2}\ge \alpha(H_v)=\alpha(F_{k-1}(G-v))\ge n_2
\] 
and 
\[
\frac{\sum_{(A,B)\in E(H_v')}((\y_v')_A-(\y_v')_B)^2}{\sum_{A\in V(H_v')}((\y_v')_A)^2}\ge \alpha(H_v')=\alpha(F_k(G-v))=n_2+1.
\]
Thus, 
    \begin{eqnarray*}
        \alpha(F_{k}(K^*_{n_1,n_2}))&=&\w^{\top}\L(F_k(K^*_{n_1,n_2}))\w\\
        &=&\sum_{(A,B)\in E(F_k(K^*_{n_1,n_2}))}(\w_A-\w_B)^2\\
        &\ge &\sum_{(A,B)\in E(H_v)}\left((\y_v)_A-(\y_v)_B\right)^2+\sum_{(A,B)\in E(H_v')}\left((\y_v')_A-(\y_v')_B\right)^2\\
        &\ge &\alpha(F_{k-1}(G-v))\sum_{A\in V(H_v)}\left((\y_v)_A\right)^2+\alpha(F_k(G-v))\sum_{A\in V(H_v')}\left((\y_v')_A\right)^2\\
        &\ge &n_2\left[\sum_{A\in V(H_v)}(\w_A)^2+\sum_{A\in V(H_v')}(\w_A)^2\right]\\
        &=&n_2. 
    \end{eqnarray*}
    For $n_1>2$ and $k>1$, we obtain the claim using a similar approach as above. 
    Together with $\alpha(F_k(K^*_{n_1,n_2}))\le n_2$, the result holds.
    \qed
\end{proof}


\section{Kite graphs}
\label{sec:kite}
The graphs in Lemma \ref{le-Shao2008} suggest the following definition. 

\begin{definition}
A {\em kite} $K(H,T)$ is a graph consisting of a {\em head} $H$ that is a general graph, 
and a {\em tail} $T$ consisting of $s\ge 2$ paths of equal length $r$ with the common end-vertex $v\in V(H)$.
\end{definition}

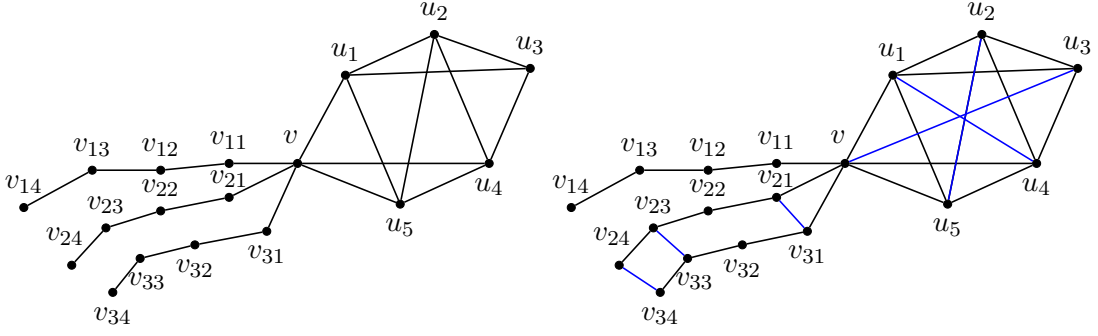
\begin{figure}
    \begin{center}
    \begin{tikzpicture}[scale=0.9,auto,swap]
    \vertex (v) at (0,1) [fill]{};
    \node () at (-0.1,1) 
    [label=above:{$v$}]{};
    \vertex (11) at (-1,1) [fill,label=above:{$v_{11}$}]{};
    \vertex (12) at (-2,0.9) [fill,label=above:{$v_{12}$}]{};
    \vertex (13) at (-3,0.9) [fill,label=above:{$v_{13}$}]{};
    \vertex (14) at (-4,0.35) [fill,label=above:{$v_{14}$}]{};
    \vertex (21) at (-1,0.5) [fill]{};
    \node () at (-1,0.35) [label=above:{$v_{21}$}]{};
    \vertex (22) at (-2,0.3) [fill,label=above:{$v_{22}$}]{};
    \vertex (23) at (-2.8,0.05) [fill,label=above:{$v_{23}$}]{};
    \vertex (24) at (-3.3,-0.5) [fill]{};
    \node () at (-3.4,-0.5) [label=above:{$v_{24}$}]{};
    \vertex (31) at (-0.45,0) [fill,label=below:{$v_{31}$}]{};
    \vertex (32) at (-1.5,-0.2) [fill,label=below:{$v_{32}$}]{};
    \vertex (33) at (-2.3,-0.4) [fill]{};
    \node () at (-2.2,-0.3) [label=below:{$v_{33}$}]{};
    \vertex (34) at (-2.7,-0.9) [fill,label=below:{$v_{34}$}]{};
    \vertex (1) at (0.7,2.3) [fill,label=above:{$u_{1}$}]{};
    \vertex (2) at (2,2.9) [fill,label=above:{$u_{2}$}]{};
    \vertex (3) at (3.4,2.4) [fill,label=above:{$u_{3}$}]{};
    \vertex (4) at (2.8,1) [fill,label=below:{$u_{4}$}]{};
    \vertex (5) at (1.5,0.4) [fill,label=below:{$u_{5}$}]{};
    \draw[line width=0.6pt](v)--(11);
    \draw[line width=0.6pt](11)--(12);
    \draw[line width=0.6pt](12)--(13);
    \draw[line width=0.6pt](13)--(14);
    \draw[line width=0.6pt](v)--(21);
    \draw[line width=0.6pt](21)--(22);
    \draw[line width=0.6pt](22)--(23);
    \draw[line width=0.6pt](23)--(24);
    \draw[line width=0.6pt](v)--(31);
    \draw[line width=0.6pt](31)--(32);
    \draw[line width=0.6pt](32)--(33);
    \draw[line width=0.6pt](33)--(34);
    \draw[line width=0.6pt](v)--(1);
    \draw[line width=0.6pt](1)--(2);
    \draw[line width=0.6pt](2)--(3);
    \draw[line width=0.6pt](3)--(4);
    \draw[line width=0.6pt](4)--(5);
    \draw[line width=0.6pt](v)--(5);
    \draw[line width=0.6pt](1)--(5);
    \draw[line width=0.6pt](1)--(3);
    \draw[line width=0.6pt](2)--(4);
    \draw[line width=0.6pt](2)--(5);
    \draw[line width=0.6pt](v)--(4);
    \vertex (v') at (8,1) [fill]{};
    \node () at (7.9,1)
    [label=above:{$v$}]{};
    \vertex (11') at (7,1) [fill,label=above:{$v_{11}$}]{};
    \vertex (12') at (6,0.9) [fill,label=above:{$v_{12}$}]{};
    \vertex (13') at (5,0.9) [fill,label=above:{$v_{13}$}]{};
    \vertex (14') at (4,0.35) [fill,label=above:{$v_{14}$}]{};
    \vertex (21') at (7,0.5) [fill]{};
    \node () at (7,0.35) [label=above:{$v_{21}$}]{};
    \vertex (22') at (6,0.3) [fill,label=above:{$v_{22}$}]{};
    \vertex (23') at (5.2,0.05) [fill,label=above:{$v_{23}$}]{};
    \vertex (24') at (4.7,-0.5) [fill]{};
    \node () at (4.6,-0.5) [label=above:{$v_{24}$}]{};
    \vertex (31') at (7.45,0) [fill,label=below:{$v_{31}$}]{};
    \vertex (32') at (6.5,-0.2) [fill,label=below:{$v_{32}$}]{};
    \vertex (33') at (5.7,-0.4) [fill]{};
    \node () at (5.8,-0.3)
    [label=below:{$v_{33}$}]{};
    \vertex (34') at (5.3,-0.9) [fill,label=below:{$v_{34}$}]{};
    \vertex (1') at (8.7,2.3) [fill,label=above:{$u_1$}]{};
    \vertex (2') at (10,2.9) [fill,label=above:{$u_2$}]{};
    \vertex (3') at (11.4,2.4) [fill,label=above:{$u_3$}]{};
    \vertex (4') at (10.8,1) [fill,label=below:{$u_{4}$}]{};
    \vertex (5') at (9.5,0.4) [fill,label=below:{$u_{5}$}]{};
    \draw[line width=0.6pt](v')--(11');
    \draw[line width=0.6pt](11')--(12');
    \draw[line width=0.6pt](12')--(13');
    \draw[line width=0.6pt](13')--(14');
    \draw[line width=0.6pt](v')--(21');
    \draw[line width=0.6pt](21')--(22');
    \draw[line width=0.6pt](22')--(23');
    \draw[line width=0.6pt](23')--(24');
    \draw[line width=0.6pt](v')--(31');
    \draw[line width=0.6pt](31')--(32');
    \draw[line width=0.6pt](32')--(33');
    \draw[line width=0.6pt](33')--(34');
    \draw[line width=0.6pt](v')--(1');
    \draw[line width=0.6pt](1')--(2');
    \draw[line width=0.6pt](2')--(3');
    \draw[line width=0.6pt](3')--(4');
    \draw[line width=0.6pt](4')--(5');
    \draw[line width=0.6pt](v')--(5');
    \draw[line width=0.6pt](1')--(5');
    \draw[line width=0.6pt](1')--(3');
    \draw[line width=0.6pt](2')--(4');
    \draw[line width=0.6pt](2')--(5');
    \draw[line width=0.6pt](v')--(4');
    \draw[blue,line width=0.6pt](21')--(31');
    \draw[blue,line width=0.6pt](23')--(33');
    \draw[blue,line width=0.6pt](24')--(34');
    \draw[blue,line width=0.6pt](3')--(v');
    \draw[blue,line width=0.6pt](2')--(5');
    \draw[blue,line width=0.6pt](4')--(1');
    \end{tikzpicture}
    \caption{A kite graph $G=K(H,T)$ and its extended graph $G^+$.}\label{fig:kite}
    \end{center}
\end{figure}

Examples of kite graphs are shown in Figures \ref{fig:Kite2}, \ref{fig:kite}, and \ref{fig:ex-th-kite}. Thus, general kite graphs appear in Lemma \ref{le-Shao2008}, whereas the graphs in Theorem \ref{th-addedge-pendant} correspond to kite graphs with tail length $r=1$.
Before proving the main theorems about kite graphs, we first need some results.
\begin{lemma}[Bapat and Pati \cite{Bapat_1998}]
\label{le-Bapat1998}
    Let $G$ be a connected graph. Let $W$ be a set of vertices of $G$ such that $G- W$ is disconnected. 
    Let $G_1$ and $G_2$ be two components of $G-W$. 
    Suppose that $\lambda_1(\L_{[V(G_1)]}(G))\le \lambda_1(\L_{[V(G_2)]}(G))$. Then, one of the two following conditions holds:
    \begin{itemize}
\item[$(i)$] 
$\alpha(G)< \lambda_1(\L_{[V(G_2)]}(G))$;
\item[$(ii)$]   
$\alpha(G)=\lambda_1(\L_{[V(G_1)]}(G))= \lambda_1(\L_{[V(G_2)]}(G))$.
    \end{itemize} 
    \end{lemma}
The following result gives a sufficient condition for having the equality
in Lemma \ref{le-Bapat1998}$(ii)$.
\begin{lemma}\label{le-suffcondition}
    Let $G$ be a connected graph with a cut vertex $v$. 
    Let $G_1,G_2,\ldots,G_r$ be $r$ components of $G-v$ with $\lambda_1(\L_{[V(G_1)]}(G))\le \lambda_1(\L_{[V(G_2)]}(G))\le \cdots \le\lambda_1(\L_{[V(G_r)]}(G))$. 
    If $\lambda_1(\L_{[V(G_1)]}(G))=\lambda_1(\L_{[V(G_2)]}(G))$, then $\alpha(G)=\lambda_1(\L_{[V(G_1)]}(G))=\lambda_1(\L_{[V(G_2)]}(G))$. 
\end{lemma}
\begin{proof}
It follows from Lemma \ref{le-principal_submatrix} that 
\begin{equation}
\label{eq:ineq-lambdas}
\lambda_1(\L(G))\le \lambda_1(\L_{[V(G-v)]}(G))\le \lambda_2(\L(G))\le \lambda_2(\L_{[V(G-v)]}(G)).
\end{equation}
Since $\lambda_1(\L_{[V(G_1)]}(G))=\lambda_1(\L_{[V(G_2)]}(G))$, we have $\lambda_1(\L_{[V(G-v)]}(G))= \lambda_2(\L_{[V(G-v)]}(G))$. 
This is because $\L_{[V(G-v)]}(G)$ is a block diagonal matrix with blocks $\L_{[V(G_i)]}(G)$ for $i=1,\ldots,r$ and, hence, the spectrum of the former is the union of the spectra of the latter.
Then, \eqref{eq:ineq-lambdas} gives $\alpha(G)=\lambda_2(\L(G))=\lambda_1(\L_{[V(G_1)]}(G))=\lambda_1(\L_{[V(G_2)]}(G))$. \qed
\end{proof}
We next give a sufficient and necessary condition for $\alpha (G)=2+2\cos(\frac{2r\pi}{2r+1})$ if $G$ is a kite graph.
\begin{lemma}
\label{le-kite-iff}
    Let $G=K(H,T)$ be a kite graph with a tail of length $r$ formed with $s(\ge 2)$ paths. Then,  $\alpha (G)=2+2\cos(\frac{2r\pi}{2r+1})$ if and only if $\lambda_1(\L_{[V(H-v)]}(H))\ge 2+2\cos(\frac{2r\pi}{2r+1}).$
\end{lemma}
\begin{proof}
    Note that \[\L_{[V(G-v)]}(G)=\diag\left(\L_{[V(H-v)]}(H),\L_{[V_{r}]}(P_{r+1}),\ldots,\L_{[V_{r}]}(P_{r+1})\right),\]  
    and $\lambda_1(\L_{[V_{r}]}(P_{r+1}))=2+2\cos(\frac{2r\pi}{2r+1})$. 
    If $\lambda_1(\L_{[V(H-v)]}(H))\ge \lambda_1(\L_{[V_{r}]}(P_{r+1}))$, then it follows from Lemma \ref{le-suffcondition} that $\alpha(G)=\lambda_1(\L_{[V_{r}]}(P_{r+1}))=2+2\cos(\frac{2r\pi}{2r+1})$.
    
Conversely, suppose, by way of contradiction, that $\lambda_1(\L_{[V(H-v)]}(H))<\lambda_1(\L_{[V_{r}]}(P_{r+1}))$. Then, it follows from Lemma 
\ref{le-Bapat1998} that $\alpha(G)<\lambda_1(\L_{[V_{r}]}(P_{r+1}))=2+2\cos(\frac{2r\pi}{2r+1})$, a contradiction.
Hence, we obtain $\lambda_1(\L_{[V(H-v)]}(H))\ge 2+2\cos(\frac{2r\pi}{2r+1})$.\qed
\end{proof}


\begin{proposition}
\label{prop:LvsS}
 Let $G$ be a kite graph $K(H,T)$ with head $H$ and tail $T$ with $s(\ge 2)$ paths of length $r$, with vertices $v_{i1},v_{i2},\ldots,v_{ir}$ for $i=1,\ldots,s$. Let $G^+$ be a graph obtained from $G$ by adding some edges among the vertices of each set $U_j=\{v_{2j},v_{3j},\ldots,v_{sj}\}$ for each given $j\in \{1,\ldots,r\}$. Let $\L$ be the Laplacian matrix of $G$. Let $\SS=(s_{ij})$ be the matrix, indexed also by the vertices of $G$, with entries
 $$
 s_{uv}=\left\{
 \begin{array}{lll}
 1 & \mbox{if } u=v\in V(H)\cup\{v_{11},v_{12}, \ldots,v_{1r}\},\\
\frac{1}{s-1} & \mbox{if } u,v\in U_j\quad \mbox{for every } j=1,\ldots,r,\\
 0 & \mbox{elsewhere.}
 \end{array}
 \right.
 $$
 Then, the following statements hold. 
 \begin{itemize}
\item[$(i)$]
The matrices $\L$ and $\SS$ commute: $\L\SS=\SS\L$.
\item[$(ii)$]
For every eigenvalue $\lambda$ of $G$, there exists a corresponding eigenvector  $\x$ of $G$ such that $\x_{v_{2j}}=\x_{v_{3j}}=\cdots=\x_{v_{sj}}$ for $1\le j\le r$.
\item[$(iii)$]
Every different eigenvalue of $G$ appears also as an eigenvalue of 
$G^{+}$ (not necessarily including multiplicities). In particular, $\alpha(G^{+})=\alpha(G)$.
\end{itemize}
 \end{proposition}
 \begin{proof}
 $(i)$ For simplicity, assume that $H=\{v\}$ (if not, the proof is basically the same). Then, the block matrices $\L$ and $\SS$ are
 $$
  \L=
  {\scriptsize{\left(
  \begin{array}{cccccc}
  d_v & -\j^{\top} & \vec0 & \vec0 & \cdots & \vec0\\
  -\j & 2\I & -\I & \O & \cdots & \O\\
  \vec0 & -\I & \ddots & \ddots & \cdots & \vdots\\
  \vec0 & \O & \ddots & \ddots & -\I & \O\\
  \vdots & \vdots & \vdots & -\I & 2\I & -\I \\
  \vec0 & \O & \cdots & \O & -\I & \I 
  \end{array}
  \right)}},\quad\mbox{and}\quad
  \SS=
  \diag(1,\SS_{11},\ldots,\SS_{rr})
  $$
  where  $\SS_{11}=\cdots=\SS_{rr}=\diag\left(1,\frac{1}{s-1}\J\right)$. Then, the computation of $\L\SS$ and $\SS\L$ shows that the commutativity comes from $\j^{\top}\SS_{11}=\j^{\top}$ and $\SS_{11}\j=\j$.\\
 $(ii)$  Let $\y$ be a $\lambda$-eigenvector of $G$ so that $\L\y=\lambda\y$. First, notice that the vector $\x=\SS\y$ satisfy the condition since
    \begin{align}
    \x_{v} &=\y_{v},\quad\text{ if } u\in V(H)\cup\{v_{11},v_{12},\ldots,v_{1r}\};
     \label{eq-kite-vec1}\\
    \x_{v_{2j}} &=\x_{v_{3j}}=\cdots=\x_{v_{sj}}=\frac{\y_{v_{2j}}+\cdots+\y_{v_{sj}}}{s-1},\quad \text{ if } j=1,\ldots, r.
    \label{eq-kite-vec2}
    \end{align}
Moreover, since $\y\neq \vec0$ and all the paths of the tail are isomorphic, there is no loss of generality if we assume that $\y_u\neq 0$ for some $u$ in \eqref{eq-kite-vec1}. Then,
$\x\neq \vec0$, and
  $$
  \L\x=\L\SS\y=\SS\L\y=\lambda\SS\y=\lambda\x,
  $$
  so that $\x$ is also a $\lambda$-eigenvector of $G$.\\
  $(iii)$ As happens for the algebraic connectivity in Lemma \ref{le-adding_edge_iff}, if a $\lambda$-eigenvector $\x$ of $G$ has equal components $\x_u=\x_v$, as in \eqref{eq-kite-vec2}, the graph $G^{+}$, obtained by adding the edge $uv$,
  has the same eigenvector $\x$ and eigenvalue $\lambda$. \qed
 \end{proof}

\begin{example}
The kite graph $G$ of Figure \ref{fig:Kite2} (left) 
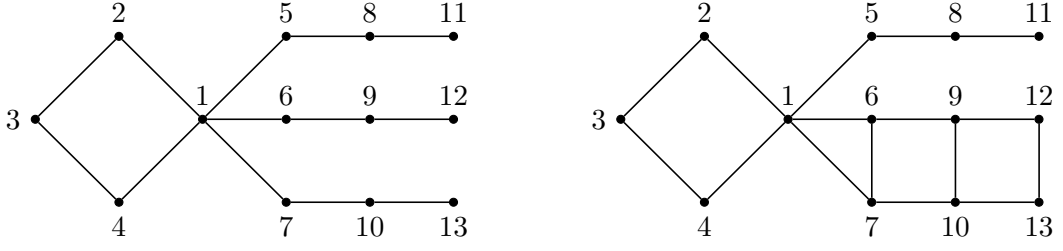
\begin{figure}
    \begin{center}
    \begin{tikzpicture}[scale=1.1,auto,swap]
    \vertex (1) at (0,1) [fill,label=above:{$1$}]{};
    \vertex (2) at (-1,2) [fill,label=above:{$2$}]{};
    \vertex (3) at (-2,1) [fill,label=left:{$3$}]{};
    \vertex (4) at (-1,0) [fill,label=below:{$4$}]{};   
    \vertex (5) at (1,2) [fill,label=above:{$5$}]{};
    \vertex (8) at (2,2) [fill,label=above:{$8$}]{};
    \vertex (11) at (3,2) [fill,label=above:{$11$}]{};
    \vertex (6) at (1,1) [fill,label=above:{$6$}]{};
    \vertex (9) at (2,1) [fill,label=above:{$9$}]{};
    \vertex (12) at (3,1) [fill,label=above:{$12$}]{};
    \vertex (7) at (1,0) [fill,label=below:{$7$}]{};
    \vertex (10) at (2,0) [fill,label=below:{$10$}]{};
    \vertex (13) at (3,0) [fill,label=below:{$13$}]{};
    \draw[line width=0.6pt](1)--(2);
    \draw[line width=0.6pt](2)--(3);
    \draw[line width=0.6pt](3)--(4);
    \draw[line width=0.6pt](1)--(4);
    \draw[line width=0.6pt](1)--(5);
    \draw[line width=0.6pt](1)--(6);
    \draw[line width=0.6pt](1)--(7);
    \draw[line width=0.6pt](5)--(8);
    \draw[line width=0.6pt](8)--(11);
    \draw[line width=0.6pt](6)--(9);
    \draw[line width=0.6pt](9)--(12);
    \draw[line width=0.6pt](7)--(10);
    \draw[line width=0.6pt](10)--(13);
    \vertex (1') at (7,1) [fill,label=above:{$1$}]{};
    \vertex (2') at (6,2) [fill,label=above:{$2$}]{};
    \vertex (3') at (5,1) [fill,label=left:{$3$}]{};
    \vertex (4') at (6,0) [fill,label=below:{$4$}]{};   
    \vertex (5') at (8,2) [fill,label=above:{$5$}]{};
    \vertex (8') at (9,2) [fill,label=above:{$8$}]{};
    \vertex (11') at (10,2) [fill,label=above:{$11$}]{};
    \vertex (6') at (8,1) [fill,label=above:{$6$}]{};
    \vertex (9') at (9,1) [fill,label=above:{$9$}]{};
    \vertex (12') at (10,1) [fill,label=above:{$12$}]{};
    \vertex (7') at (8,0) [fill,label=below:{$7$}]{};
    \vertex (10') at (9,0) [fill,label=below:{$10$}]{};
    \vertex (13') at (10,0) [fill,label=below:{$13$}]{};
    \draw[line width=0.6pt](1')--(2');
    \draw[line width=0.6pt](2')--(3');
    \draw[line width=0.6pt](3')--(4');
    \draw[line width=0.6pt](1')--(4');
    \draw[line width=0.6pt](1')--(5');
    \draw[line width=0.6pt](1')--(6');
    \draw[line width=0.6pt](1')--(7');
    \draw[line width=0.6pt](5')--(8');
    \draw[line width=0.6pt](8')--(11');
    \draw[line width=0.6pt](6')--(9');
    \draw[line width=0.6pt](9')--(12');
    \draw[line width=0.6pt](7')--(10');
    \draw[line width=0.6pt](10')--(13');
    \draw[line width=0.6pt](12')--(13');
    \draw[line width=0.6pt](6')--(7');
    \draw[line width=0.6pt](9')--(10');
    \end{tikzpicture}
    \caption{A kite graph $G=(H,T)$ with head $H=C_4$ and tail $T$ with $s=3$ paths of length $r=3$, and its extended graph $G^+$.}\label{fig:Kite2}
    \end{center}
\end{figure}
has Laplacian matrix 
$$
\L = 
\left({\scriptsize\begin{array}{ccccccccccccc}
5 & -1 & 0 & -1 & -1 & -1 & -1 & 0 & 0 & 0 & 0 & 0 & 0 
\\
 -1 & 2 & -1 & 0 & 0 & 0 & 0 & 0 & 0 & 0 & 0 & 0 & 0 
\\
 0 & -1 & 2 & -1 & 0 & 0 & 0 & 0 & 0 & 0 & 0 & 0 & 0 
\\
 -1 & 0 & -1 & 2 & 0 & 0 & 0 & 0 & 0 & 0 & 0 & 0 & 0 
\\
 -1 & 0 & 0 & 0 & 2 & 0 & 0 & -1 & 0 & 0 & 0 & 0 & 0 
\\
 -1 & 0 & 0 & 0 & 0 & 2 & 0 & 0 & -1 & 0 & 0 & 0 & 0 
\\
 -1 & 0 & 0 & 0 & 0 & 0 & 2 & 0 & 0 & -1 & 0 & 0 & 0 
\\
 0 & 0 & 0 & 0 & -1 & 0 & 0 & 2 & 0 & 0 & -1 & 0 & 0 
\\
 0 & 0 & 0 & 0 & 0 & -1 & 0 & 0 & 2 & 0 & 0 & -1 & 0 
\\
 0 & 0 & 0 & 0 & 0 & 0 & -1 & 0 & 0 & 2 & 0 & 0 & -1 
\\
 0 & 0 & 0 & 0 & 0 & 0 & 0 & -1 & 0 & 0 & 1 & 0 & 0 
\\
 0 & 0 & 0 & 0 & 0 & 0 & 0 & 0 & -1 & 0 & 0 & 1 & 0 
\\
 0 & 0 & 0 & 0 & 0 & 0 & 0 & 0 & 0 & -1 & 0 & 0 & 1 
\end{array}}\right)
$$
with algebraic connectivity $\alpha(G)=0.19086$. Then,
we can add the edges $2\sim 3$, $5\sim 6$, or $8\sim 9$ to obtain $G^+$, shown in the same figure (right), keeping the same value of $\alpha(G)$, as shown in Proposition \ref{prop:LvsS}.  Thus, a Field vector is $\y =(
0,0,0,0,-0.44504,0,0.44504,$ $- 0.80193,0,0.80193,-1,0,1)^{\top}$. If we multiply it by the matrix  $\SS$ 
$$
\SS = 
\\
\left({\scriptsize\begin{array}{ccccccccccccc}
1 & 0 & 0 & 0 & 0 & 0 & 0 & 0 & 0 & 0 & 0 & 0 & 0 
\\
 0 & 1 & 0 & 0 & 0 & 0 & 0 & 0 & 0 & 0 & 0 & 0 & 0 
\\
 0 & 0 & 1 & 0 & 0 & 0 & 0 & 0 & 0 & 0 & 0 & 0 & 0 
\\
 0 & 0 & 0 & 1 & 0 & 0 & 0 & 0 & 0 & 0 & 0 & 0 & 0 
\\
 0 & 0 & 0 & 0 & 1 & 0 & 0 & 0 & 0 & 0 & 0 & 0 & 0 
\\
 0 & 0 & 0 & 0 & 0 & \frac{1}{2} & \frac{1}{2} & 0 & 0 & 0 & 0 & 0 & 0 
\\
 0 & 0 & 0 & 0 & 0 & \frac{1}{2} & \frac{1}{2} & 0 & 0 & 0 & 0 & 0 & 0 
\\
 0 & 0 & 0 & 0 & 0 & 0 & 0 & 1 & 0 & 0 & 0 & 0 & 0 
\\
 0 & 0 & 0 & 0 & 0 & 0 & 0 & 0 & \frac{1}{2} & \frac{1}{2} & 0 & 0 & 0 
\\
 0 & 0 & 0 & 0 & 0 & 0 & 0 & 0 & \frac{1}{2} & \frac{1}{2} & 0 & 0 & 0 
\\
 0 & 0 & 0 & 0 & 0 & 0 & 0 & 0 & 0 & 0 & 1 & 0 & 0 
\\
 0 & 0 & 0 & 0 & 0 & 0 & 0 & 0 & 0 & 0 & 0 & \frac{1}{2} & \frac{1}{2} 
\\
 0 & 0 & 0 & 0 & 0 & 0 & 0 & 0 & 0 & 0 & 0 & \frac{1}{2} & \frac{1}{2} 
\end{array}}\right),
$$
then $\L\SS=\SS\L$ and 
we get
\begin{align*}
&\x=\SS\y=\\
&(0,0,0,0,-0.44504,0.22252,0.22252,-0.80193,0.40096,0.40096,-1,0.5,0.5)^{\top}.
\end{align*}
 In fact, every eigenvalue in
 $$
 \spec G=\{6.271,
3.342,
3.246^2,
2.810,
2,
1.554^2,
1.211,
0.364,
0.198^2,0\}
$$
appears (in boldface) in 
\begin{align*}
    &\spec G^+=\\
&\{{\bf 6.271},
5.246,
3.554,
{\bf 3.342},
{\bf 3.246},
{\bf 2.810},
2.198,
{\bf 2},
{\bf 1.554},
{\bf 1.211},
{\bf 0.364},
{\bf 0.198},
{\bf 0}\}.
\end{align*}
\end{example}

Now, we are ready to obtain the result on the token graph of kite graphs.
\begin{theorem}
\label{th-kite}
    Let $G=K(H,T)$ be a kite graph with a tail of length $r$ formed with $s(\ge 2)$ paths $v,v_{i1},v_{i2},\ldots,v_{ir}$ for $i=1,\ldots,s$. Let $K(H,T^+)$ be the graph obtained from $G$ by adding some edges between the vertices $v_{1j},v_{2j},\ldots,v_{sj}$, for each $j$ of $T$, and $K^+(H,T)$ be the graph obtained from $G$ by adding some edges to $H$, and between the vertices $v_{2j},\ldots,v_{sj}$ with $1\le j\le r$. Suppose that $\alpha(H)=\alpha(F_k(H))$, then the following statements hold.
\begin{itemize}
\item[$(i)$] 
If $\alpha (G)\neq2+2\cos(\frac{2r\pi}{2r+1})$, then $\alpha(F_k(K(H,T^+)))=\alpha(K(H,T^+))=\alpha(G)$. 
\item[$(ii)$] 
If $\alpha (G)=2+2\cos(\frac{2r\pi}{2r+1})$, then $\alpha(F_k(K^+(H,T)))=\alpha(K^+(H,T))=\alpha(G)=2+2\cos(\frac{2r\pi}{2r+1})$.
\end{itemize}
\end{theorem}
\begin{proof}
($i$) Note that the result by Shao, Guo, and Shan \cite{Shao_2008} in Lemma \ref{le-Shao2008} reads as $\alpha(K(H,T^+))=\alpha(G)$. 
Combining it with Theorems \ref{th-add_edge} and \ref{th-Dalfo2024}$(ii)$, and the hypothesis, we obtain $\alpha(F_k(K(H,T^+)))=\alpha(K(H,T^+))=\alpha(G)$.

($ii$) Let $G'\cong K(H^+,T)$ be the graph obtained from $G$ by adding some edges to $H$, where $|V(H)|=h$. We first prove that $\alpha(K(H^+,T))=\alpha(K(H,T))$. 
From Lemma \ref{le-kite-iff}, we get that $\lambda_1(\L_{[V(H-v)]}(H))\ge \lambda_1(\L_{[V_{r}]}(P_{r+1}))=2+2\cos(\frac{2r\pi}{2r+1})$.
Besides, observe that $G'$ is the graph obtained by adding some edges in $H$. 
Denote these additional edges by $E'$ if they are between the vertices of $H-v$, and by $E''$ if they have the end-vertex $v$. 
Let $H'$  be the graph with vertex set $V(H')=V(H)$ and edge set $E'$.
Then,  
$\L_{[V(H-v)]}(G')=\L_{[V(H-v)]}(H)+\L(H')+\I''$,
where $\I''$ is a diagonal matrix of order $h-1$, with $|E''|$ ones and $h-1-|E''|(\ge 1)$ zeros.
It follows from Lemma \ref{le-matrix_add} with $i=1$ and the hypothesis that 
\begin{align}
\lambda_1(\L_{[V(H-v)]}(G')) &\ge \lambda_1(\L_{[V(H-v)]}(H))+\lambda_1(\L(H'))+\lambda_1(\I'')\nonumber\\
 & \ge\lambda_1(\L_{[V(H-v)]}(H))\ge2+2\cos\left(\frac{2r\pi}{2r+1}\right).
\label{eq:2nd-ineq}
\end{align}
    Then, from 
    \eqref{eq:2nd-ineq} and  Lemma \ref{le-suffcondition}, we have that 
    $$
    \alpha(G)=2+2\cos\left(\frac{2r\pi}{2r+1}\right)=\alpha(G').
    $$

    Next, let $G^+\cong K^+(H,T)$ be the graph obtained from $G'$ by adding some edges between the vertices $v_{2j},\ldots,v_{sj}$ with $1\le j\le r$. It follows from Proposition \ref{prop:LvsS}$(iii)$ that $\alpha(G^+)=\alpha(G')=\alpha(G)$.  
    Together with Theorem \ref{th-Dalfo2024}$(ii)$ (since $\alpha(H)=\alpha(F_k(H))$ implies $\alpha(G)=\alpha(F_k(G))$) and Theorem \ref{th-add_edge}, we obtain $\alpha(F_k(K^+(H,T)))=\alpha(K^+(H,T))=\alpha(G)=2+2\cos\left(\frac{2r\pi}{2r+1}\right)$.
    \qed
\end{proof}

\begin{example}
    Figure \ref{fig:ex-th-kite} is an example of 
    \eqref{eq:2nd-ineq} with $E'=\{36,46\}$ and $E''=\{13,14\}$. Note that $\L_{[V(H-v)]}(G^+)=\L_{[V(H-v)]}(H)+\L(H')+\I''$ with 
    \begin{eqnarray*}
        \L_{[V(H-v)]}(H)
        &=&{\scriptsize{\left({\begin{array}{ccccc}
        3&-1&0&-1&0\\
        -1&2&-1&0&0\\
        0&-1&2&-1&0\\
        -1&0&-1&3&-1\\
        0&0&0&-1&2
        \end{array}}
        \right)}},\\
        \L(H')&=&{\scriptsize{\left({\begin{array}{ccccc}
        0&0&0&0&0\\
        0&1&0&0&-1\\
        0&0&1&0&-1\\
        0&0&0&0&0\\
        0&-1&-1&0&2
        \end{array}}
        \right)}},
        ~\text{and } \I''={\scriptsize{\left({\begin{array}{ccccc}
        0&0&0&0&0\\
        0&1&0&0&0\\
        0&0&1&0&0\\
        0&0&0&0&0\\
        0&0&0&0&0
        \end{array}}
        \right).}}
    \end{eqnarray*}
    Moreover, $\lambda_1(\L_{[V(H-v)]}(G^+))=0.747>0.284=\lambda_1(\L_{[V(H-v)]}(H))$.
\end{example}
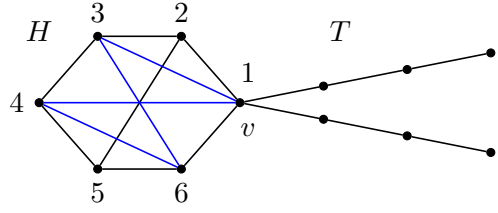
\begin{figure}
    \begin{center}
    \begin{tikzpicture}[scale=1.1,auto,swap]
    \vertex (1) at (0,1) [fill]{};
    \node () at (0.1,1) [label=above:{$1$}]{};
    \node () at (0.1,1) [label=below:{$v$}]{};
    \vertex (01) at (1,1.2) [fill]{};
    \vertex (02) at (2,1.4) [fill]{};
    \vertex (03) at (3,1.6) [fill]{};
    \vertex (04) at (1,0.8) [fill]{};
    \vertex (05) at (2,0.6) [fill]{};
    \vertex (06) at (3,0.4) [fill]{};
    \vertex (2) at (-0.7,1.8) [fill,label=above:{$2$}]{};
    \vertex (3) at (-1.7,1.8) [fill,label=above:{$3$}]{};
    \vertex (4) at (-2.4,1) [fill,label=left:{$4$}]{};
    \vertex (5) at (-1.7,0.2) [fill,label=below:{$5$}]{};
    \vertex (6) at (-0.7,0.2) [fill,label=below:{$6$}]{};   
    \node () at (-2.4,1.5) [label=above:{$H$}]{};
    \node () at (1.2,1.5) [label=above:{$T$}]{};
    \draw[line width=0.6pt](1)--(01);
    \draw[line width=0.6pt](01)--(02);
    \draw[line width=0.6pt](02)--(03);
    \draw[line width=0.6pt](1)--(04);
    \draw[line width=0.6pt](04)--(05);
    \draw[line width=0.6pt](05)--(06);
    \draw[line width=0.6pt](1)--(2);
    \draw[line width=0.6pt](2)--(3);
    \draw[line width=0.6pt](3)--(4);
    \draw[line width=0.6pt](4)--(5);
    \draw[line width=0.6pt](5)--(6);
    \draw[line width=0.6pt](1)--(6);
    \draw[line width=0.6pt](2)--(5);
    \draw[blue,line width=0.6pt](1)--(3);
    \draw[blue,line width=0.6pt](1)--(4);
    \draw[blue,line width=0.6pt](3)--(6);
    \draw[blue,line width=0.6pt](4)--(6);
    \end{tikzpicture}
    \caption{A example for Theorem \ref{th-kite}.}\label{fig:ex-th-kite}
    \end{center}
\end{figure}

This theorem and the following ones in this section can be generalized for superkite graphs. A superkite graph $K(H,{\cal T})$ has head $H$ and supertail ${\cal T}=Q^s$ formed by $s(\ge 2)$ trees isomorphic to $Q$ with the same root $v\in V(H)$. Then, basically, the same proof of Theorem \ref{th-kite} leads to the following result.
\begin{theorem}
    \label{th-superkite}
    Let $G$ be a superkite graph $K(H,{\cal T})$ with tail ${\cal T}=Q^s$, such that $H$ is not a complete graph with $h=|V(H)|$, with $\alpha(H)=\alpha(F_k(H))$, and 
    $$
    \lambda_1(\L_{[V(H-v)]}(H))\ge \lambda_1(\L_{[V(Q-v)]}(Q)).
    $$
    Let $G^+$ be the graph obtained from $G$ by adding some edges to $H$. Then, $\alpha(G^+)=\alpha(F_k(G^+))$.
\end{theorem}
The following two corollaries deal with the cases in which the head of a kite graph is a cycle and a complete bipartite graph, respectively.
\begin{corollary}
\label{th-kite-cycle}
    Let $G$ be a kite graph $K(H,T)$ with tail length $r$, such that $H\cong C_h$ is a cycle with order $h$ $(\mbox{for } h\le 2r+1)$ and $\alpha(H)=\alpha(F_k(H))$. 
    Let $G^+$ be the graph obtained from $G$ by adding some edges to $H$, and between the vertices $v_{2j},\ldots,v_{sj}$ with $1\le j\le r$. Then, 
    $\alpha(G^+)=\alpha(F_k(G^+))$. 
\end{corollary}
\begin{proof}
It is known that the characteristic polynomials of  $\L_{[V']}(C_{h})$, where $V'=V(C_h-v)$, and $\L(P_{h})$ are related as $x\cdot\Phi(\L_{[V']}(C_{h}))=\Phi(\L(P_{h}))$ (see Guo \cite{Guo_2008}).
Hence, $\lambda_i(\L_{[V']}(C_{h}))=\lambda_{i+1}(\L(P_{h}))$ for $i=1,\ldots,h-1$. Then, since the Laplacian eigenvalues of the path $P_{\ell}$ of length $\ell$ are $\tau_k=2-2\cos(\frac{k\pi}{\ell})$ for $k=0,1,\ldots,\ell$, we have, with $H=C_h$,
\begin{align*}
\lambda_1(\L_{[V(H-v)]}(H))=\lambda_2(\L(P_h))&\ge\lambda_2(\L(P_{2r+1}))=2-2\cos\left(\frac{\pi}{2r+1}\right)\\
 &=2+2\cos\left(\frac{2r\pi}{2r+1}\right)=
\lambda_1(\L_{[V_r]}(P_{r+1})).
    \end{align*}
    From Lemma \ref{le-kite-iff} and Theorem \ref{th-kite}, the result holds.     \qed
\end{proof}

\begin{corollary}
\label{th-kite-bipartite}
   Let $G$ be a kite graph $K(H,T)$ with root vertex $v$,  tail length $r$, and head $H\cong K_{h_1,h_2}$ being a complete bipartite graph with vertex set $V=V_1\cup V_2$, on $h=h_1+h_2$ vertices, where $h_i=|V_i|$.
   Let $G^+$ be the graph obtained from $G$ by adding some edges to $H$, and between the vertices $v_{2j},\ldots,v_{sj}$ with $1\le j\le r$. 
   If  $v\in V_i$ for some $i\in\{1,2\}$, and 
   $$
   \frac{1}{2}\left(h-\sqrt{h^2-4h_i}\right)\ge 2+2\cos\left(\frac{2r\pi}{2r+1}\right),
   $$
   then  
   $\alpha(G^+)=\alpha(F_k(G^+))$.
\end{corollary}
\begin{proof}
    It follows from Theorem \ref{th-Dalfo2024}($iii$) that $\alpha(F_k(K_{h_1,h_2}))=\alpha(K_{h_1,h_2})$.  
    Moreover, it can be shown that, if $v\in V_i$,  then $\lambda_1(\L_{V(H-v)}(H))=\frac{1}{2}\left(h-\sqrt{h^2-4h_i}\right)$ for $i=1,2$.
    Then, we have $\lambda_1(\L_{V(H-v)}(H))\ge 2+2\cos(\frac{2r\pi}{2r+1})=\lambda_1(\L_{[V_{r}]}(P_{r+1}))$.
    Therefore, the result holds by Lemma \ref{le-kite-iff} and Theorem \ref{th-kite}.
    \qed
\end{proof}

\begin{figure}
    \begin{center}
    \begin{tikzpicture}[scale=1.1,auto,swap]
        \vertex (v) at (1,1) [fill]{};
    \vertex (11) at (0,1) [fill]{};
    \vertex (12) at (-1,0.9) [fill]{};
    \vertex (13) at (-2,0.9) [fill]{};
    \vertex (21) at (0,0.5) [fill]{};
    \vertex (22) at (-1,0.3) [fill]{};
    \vertex (23) at (-1.8,0.05) [fill]{};
    \vertex (31) at (0.45,0) [fill]{};
    \vertex (32) at (-0.5,-0.2) [fill]{};
    \vertex (33) at (-1.3,-0.4) [fill]{};
    \vertex (1) at (1.8,1.9) [fill]{};
    \vertex (2) at (3,1.5) [fill]{};
    \vertex (3) at (3,0.5) [fill]{};
    \vertex (4) at (1.8,0.1) [fill]{};
    \draw[line width=0.6pt](v)--(11);
    \draw[line width=0.6pt](11)--(12);
    \draw[line width=0.6pt](12)--(13);
    \draw[line width=0.6pt](v)--(21);
    \draw[line width=0.6pt](21)--(22);
    \draw[line width=0.6pt](22)--(23);
    \draw[line width=0.6pt](v)--(31);
    \draw[line width=0.6pt](31)--(32);
    \draw[line width=0.6pt](32)--(33);
    \draw[line width=0.6pt](v)--(1);
    \draw[line width=0.6pt](1)--(2);
    \draw[line width=0.6pt](2)--(3);
    \draw[line width=0.6pt](3)--(4);
    \draw[line width=0.6pt](4)--(v);
    \vertex (2v) at (7,1) [fill]{};
    \vertex (211) at (6,1) [fill]{};
    \vertex (212) at (5,0.9) [fill]{};
    \vertex (213) at (4,0.9) [fill]{};
    \vertex (221) at (6,0.5) [fill]{};
    \vertex (222) at (5,0.3) [fill]{};
    \vertex (223) at (4.2,0.05) [fill]{};
    \vertex (231) at (6.45,0) [fill]{};
    \vertex (232) at (5.5,-0.2) [fill]{};
    \vertex (233) at (4.7,-0.4) [fill]{};
    \vertex (21) at (7.35,1.9) [fill]{};
    \vertex (22) at (8.4,1.9) [fill]{};
    \vertex (23) at (8.2,1) [fill]{};
    \vertex (24) at (7.9,0.2) [fill]{};
    \draw[line width=0.6pt](2v)--(211);
    \draw[line width=0.6pt](211)--(212);
    \draw[line width=0.6pt](212)--(213);
    \draw[line width=0.6pt](2v)--(221);
    \draw[line width=0.6pt](221)--(222);
    \draw[line width=0.6pt](222)--(223);
    \draw[line width=0.6pt](2v)--(231);
    \draw[line width=0.6pt](231)--(232);
    \draw[line width=0.6pt](232)--(233);
    \draw[line width=0.6pt](2v)--(23);
    \draw[line width=0.6pt](2v)--(22);
    \draw[line width=0.6pt](21)--(22);
    \draw[line width=0.6pt](21)--(23);
    \draw[line width=0.6pt](21)--(24);
    \draw[line width=0.6pt](24)--(2v);
    \end{tikzpicture}
    \caption{The example in Corollary \ref{th-kite-cycle} (left) and Corollary \ref{th-kite-bipartite} (right) with $r=3$ and $s=3$.}\label{fig:kite_example}
    \end{center}
\end{figure}
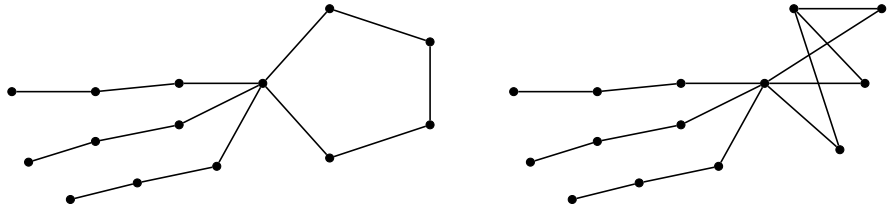

\section{Graphs with a cut clique}
\label{sec:cut-clique}
The notion of cut clique generalizes the idea of cut vertex.
Thus, a clique $K_r$ of a graph $G$ is a {\em cut clique} if the removal of (the vertices and edges of) $K_r$ disconnects $G$ in two or more components.
Then, the following result is a natural generalization of Kirland's result \cite{Kirkland_2000} for cut vertices.
\begin{proposition}
\label{prop-cut-clique}
For a connected graph $G$ with $n$ vertices and a cut clique $K_r$, we have $\alpha(G)\le r$ and the equality holds if and only if 
$d_G(K_r)=r(n-r)$.    
\end{proposition}
\begin{proof}
    It follows from Fiedler \cite[\S 3.12]{Fiedler_1973} that 
    $$
    \alpha(G)\le \min\{\alpha(G-K_r)+|V(K_r)|,\alpha(K_r)+|V(G-K_r)|\}.
    $$
    Then, since $G-K_r$ is disconnected, we have $\alpha(G)\le r$.

We next show that, if $d_G(K_r)=r(n-r)$, then $\alpha(G)=r$ .
    Let $G_1,\ldots,G_s$ be all components of $G-K_r$. 
    Notice that $\L(G)=\L_1+\L_2$, where $\L_1$ is the Laplacian matrix of the graph with vertex set $V(G)$ and edge set $E(G_1\cup\cdots\cup G_s)$, and $\L_2$ is the Laplacian matrix of $G$ removing all the edges of $G_1\cup\cdots\cup G_s$. That is,
    $$\L(G)={\scriptsize{\left({\begin{array}{ccccccc}
        \L(G_1)&&&0&\ldots&0\\
        &\ddots&&&\ddots&\\
        &&\L(G_s)&0&\ldots&0\\
        0&\ldots&0&0&\ldots&0\\
        &\ddots&& &\ddots&\\
        0&\ldots&0& 0& \ldots&0
        \end{array}}
        \right)}}+{\scriptsize{\left({\begin{array}{ccccccc}
        r\I& & &-1&\ldots&-1\\
        &\ddots & &&\ddots&\\
        &&r\I&-1&\ldots&-1\\
        -1&\ldots&-1&n-1&\ldots&-1\\
        &\ddots&& &\ddots&\\
        -1&\ldots&-1& -1& \ldots&n-1
        \end{array}}
        \right)}}.
        $$
    Since the spectrum of $\L_2$ is $\{0^1,r^{n-r-1},n^r\}$, we have that $\lambda_2(\L_2)=r$. 
    It follows from Lemma \ref{le-matrix_add} with $i=2$ that 
        $\alpha(G)=\lambda_2(\L(G))\ge\lambda_1(\L_1)+\lambda_2(\L_2)=r.$
    Together with $\alpha(G)\le r$, we conclude that $\alpha(G)=r$. 

    Otherwise, if $d_{G}(K_r)\ne r(n-r)$, then there exists $u\in V(K_r)$ and $v\in V(G_1\cup\cdots\cup G_s)$ such that $uv\notin E(G)$.
    Denote by $n_1$ and $n_2$ the order of $G_1\cup\cdots\cup G_{s-1}$ and $G_s$, respectively.
    Let $G^*$ be a graph obtained from $G$ by adding edges such that $G^*[V(G_1\cup\cdots\cup G_{s-1})]$ and $G^*[V(G_{s})]$ are complete graphs, and $d_{G^*}(K_r)=r(n-r)$.
    Note that $$\L(G^*)=\left({\begin{array}{ccc}
        (n_1+r)\I_{n_1}-\J_{n_1}&\O&-\J_{n_1\times r}\\
        \O&(n_2+r)\I_{n_2}-\J_{n_2}&-\J_{n_2\times r}\\
        -\J_{r\times n_1}&-\J_{r\times n_2}&n\I_{r}-\J_{r}
        \end{array}}
        \right).$$
    It can be shown that $\Phi(\L(G^*))=x(x-r)(x-r-n_1)^{n_1-1}
    (x-r-n_2)^{n_2-1}(x-n)^r$. Thus, 
    the spectrum of $\L(G^*)$ is $\{0^1,r^1,(r+n_1)^{n_1-1},(r+n_2)^{n_2-1},n^r\}$.
    Moreover, it is easily checked that 
    \[
    \x=\Big(\underbrace{1,\ldots,1}_{n_1},\underbrace{-\frac{n_1}{n_2},\ldots,-\frac{n_1}{n_2}}_{n_2},\underbrace{0,\ldots,0}_r\Big)^{\top}
    \]
 is a Fiedler vector of $G^*$. Furthermore, since $r$ has multiplicity 1,  any eigenvector $\y$ of $G^*$ with  eigenvalue $\alpha(G^*)$ must be a multiple of $\x$.
    It follows from Lemma \ref{le-interlacing} that $\alpha(G^*)\ge\alpha(G^*-uv)$. 
    If the equality holds and $\w$ is a Fiedler vector of $G^*-uv$, we have $\w_u=\w_v$ from Lemma \ref{le-adding_edge_iff}.
    Moreover, $\w$ also is a Fiedler vector of $G^*$ with $\w_u=\w_v$, a contradiction. 
    Hence, the inequality must be strict. 
    Finally, since $G^*-uv$ can be obtained by adding edges in $G$, we conclude that $\alpha(G)\le\alpha(G^*-uv)<\alpha(G^*)=r.$
\qed
\end{proof}
\begin{proposition}
\label{prop-subgraph}
    Let $G_1$ and $G_2$ be two edge-disjoint graphs with the same $n$ vertices. 
    Let $G=G_1\cup G_2$. Then, $F_k(G)=F_k(G_1)\cup F_k(G_2)$. 
\end{proposition}
\begin{proof}
    It suffices to consider that there is exactly one edge $uv$ in $G_2$, that is, $G_1=G-uv$.  
    The result holds since $F_k(G-uv)$ is a spanning subgraph of $F_k(G)$ with 
   \begin{align*}
    E(F_k(G))\backslash E(F_k(G-uv)) &=\{A_{r}A_{s} : A_{r}=\{u,u_1,\ldots,u_{k-1}\}, A_{s}=\{v,u_1,\ldots,u_{k-1}\}\}\\
    & =E(F_k(G_2)). \mbox{\hskip 7.8cm  \qed} 
    \end{align*}
    \end{proof}

\begin{theorem}
\label{th-clique}
Let $G$ be a graph on $n$ vertices, with a cut clique $K_r$ such that $G-K_r$ results in some disjoint graphs $G_1,\ldots,G_s$, for $s\ge 2$, and 
 degree $d_G(K_r)=r(n-r)$ (so that there are {\bf all} the edges between the vertices of $K_r$ and the vertices of each $G_i$).  Then, $\alpha(G)=\alpha(F_k(G))=r$.
\end{theorem}

\begin{proof}
    By Theorem \ref{th-Dalfo2021} and Proposition \ref{prop-cut-clique}, we get
    \begin{equation}
        \alpha(F_k(G))\le \alpha(G)=r.\label{eq-th-clique_1}
    \end{equation}
    Let $H_1$ be a graph with vertex set $V(G)$ and edge set $E(H_1)=\{uv:uv\in G_1\cup\cdots\cup G_s\}$, and let $H_2$ be a graph obtained from $G$ by removing all the edges of $E(H_1)$. Note that $H_1$ and $H_2$ are the graphs with Laplacian matrices $\L_1$ and $\L_2$, respectively, given in the proof of Proposition \ref{prop-cut-clique}, and $H_2$ is the graph obtained from a complete bipartite graph by adding all edges in the partition with order $r$.
    By Theorems \ref{th-bip-addedge1}, and \ref{th-bip-addedge2}, we obtain $\alpha(F_k(H_2))=r$. 
    From Proposition \ref{prop-subgraph}, we have $\L(F_k(G))=\L(F_k(H_1))+\L(F_k(H_2))$. 
    Then, it follows from Lemma \ref{le-matrix_add} with $i=2$ that
    \begin{equation*}
        \alpha(F_k(G))=\lambda_2(F_k(G))\ge \lambda_1(F_k(H_1))+\lambda_2(F_k(H_2))=\alpha(F_k(H_2))=r.
    \end{equation*}
    Combining it with (\ref{eq-th-clique_1}), this completes the proof.\qed
\end{proof}

In particular, with $r=1$, the vertex of $K_1$ is a cut vertex of degree $n-1$, and we reobtain a result of Barik and Verma \cite{bv24} (see Theorem \ref{theBarik} in the Preliminaries).

\end{document}